\documentclass[12pt]{article}

\usepackage{amsmath,amssymb,latexsym}

\newcommand{\SC}{\scriptstyle}

\newcommand{\bL}{\mathbf{L}}
\newcommand{\NN}{\mathbb{N}}
\newcommand{\RR}{\mathbb{R}}
\newcommand{\ZZ}{\mathbb{Z}}

\newcommand{\ee}{\mathrm{e}}

\newcommand{\Var}{\operatorname{Var}}
\newcommand{\EE}{\operatorname{\mathsf{E}}}
\newcommand{\PP}{\operatorname{\mathsf{P}}}
\newcommand{\sPP}{{\SC\PP}}

\newcommand{\hbeta}{\widehat{\beta}}
\newcommand{\hvarrho}{\widehat{\varrho}}
\newcommand{\hb}{\widehat{b}}

\newcommand{\hbb}{\widehat{\mathbf{b}}}
\newcommand{\ty}{\widetilde{y}}

\newcommand{\cD}{\mathcal{D}}

\newcommand{\cF}{\mathcal{F}}
\newcommand{\cG}{\mathcal{G}}

\newcommand{\cL}{\mathcal{L}}

\newcommand{\cN}{\mathcal{N}}

\newcommand{\scD}{{\SC\cD}}

\newcommand{\distr}{\stackrel{\scD}{\longrightarrow}}
\newcommand{\stoch}{\stackrel{\sPP}{\longrightarrow}}

\renewcommand{\leq}{\leqslant}
\renewcommand{\geq}{\geqslant}

\newcommand{\vare}{\varepsilon}

\newcommand{\proofend}{\hfill$\square$}

\newtheorem{theorem}{Theorem}
\newtheorem{lemma}{Lemma}

\newtheorem{corollary}{Corollary}
\newtheorem{remark}{Remark}

\renewcommand{\thetheorem}{\arabic{section}.\arabic{theorem}}
\renewcommand{\thelemma}{\arabic{section}.\arabic{lemma}}
\renewcommand{\theproposition}{\arabic{section}.\arabic{proposition}}
\renewcommand{\thedefinition}{\arabic{section}.\arabic{definition}}
\renewcommand{\thecorollary}{\arabic{section}.\arabic{corollary}}
\renewcommand{\theremark}{\arabic{section}.\arabic{remark}}

\title{Joint ML estimation of all parameters in a discrete time random field
       HJM type interest rate
       model}

\author{\sc J\'ozsef G\'all$^{\star}$, Gyula Pap$^{\star\star}$ and
        Martien van Zuijlen$^{\star\star\star}$ \\ 
 {\small\em  $^{\star}$Faculty of Economics, University of Debrecen, Hungary}, \\ 
 {\small jozsef.gall@econ.unideb.hu} \\ 
 {\small\em $^{\star\star}$Bolyai Institute, University of Szeged, Hungary} \\ 
 {\small papgy@math.u-szeged.hu} \\ 
 {\small\em$^{\star\star\star}$Institute for Mathematics, Astrophysics and Particle
            Physics}, \\ 
 {\small\em Radboud University, Nijmegen, The Netherlands} \\ 
 {\small M.vanZuijlen@science.ru.nl}
       }

\date{}

\begin{document}

\maketitle

\begin{abstract}  
We consider discrete time Heath--Jarrow--Morton type interest rate models,
 where the interest rate curves are driven by a geometric spatial
 autoregression field. 
Strong consistency and asymptotic normality of the maximum likelihood
 estimators of the parameters are proved for stable no-arbitrage models
 containing a general stochastic discounting factor, 
 where explicit form of the ML estimators is not available 
 given a non-i.i.d. sample. 
The results form the basis of further statistical problems in 
 such models. 

\end{abstract}

\section{Introduction}
\setcounter{theorem}{0}
\setcounter{lemma}{0}
\setcounter{proposition}{0}
\setcounter{definition}{0}
\setcounter{remark}{0}
\setcounter{corollary}{0}

Our aim in the present paper is to consider some statistical questions arising
 in a Heath--Jarrow--Morton (HJM) type interest rate model proposed by
 G\'all, Pap and Zuijlen \cite{GPZ}. 
Such models are useful not only for describing the structure of interest rates 
 but also for describing bond price structures in the market. 
We focus on asymptotic properties of the joint maximum likelihood 
 estimators (MLE) of the parameters of the model, where the 
 non-i.i.d. sample and the lack of an explicit form of the estimators 
 make the derivation of the results difficult. 
These results give the basis of further statistical problems, such as 
 hypothesis tests, interval estimations or model selection tools. 
 
In the following we specify the model. For $\ZZ_+$ being the sets of non-negative
 integers,
let $f_{k, \ell}$ denote the forward interest rate at time $k \in \ZZ_+$ with time
 to maturity date $\ell \in \ZZ_+$. 
Hence it is the interest rate  for the future time
 period $[k+\ell, k+\ell+1)$. 

The forward rate dynamics is supposed to be given by the (stochastic)
 difference equation
 \[
   f_{k + 1, \ell}
   = f_{k, \ell} + \alpha_{k, \ell}
     + \beta \Delta_1 S_{k, \ell} , \qquad k, \ell \in \ZZ_+ ,
 \]
 where the initial values $( f_{0, \ell} )_{ \ell \in \ZZ_+}$ are given real numbers,
 $\beta \in \RR$ denotes the volatility and 
 $\Delta_1 S_{k, \ell} := S_{k+1, \, \ell} -S_{k, \ell}$, where
 $( S_{k, \ell} )_{k, \ell \in \ZZ_+}$ is a doubly geometric spatial autoregressive
 process given by
 \[
   \begin{cases}
    S_{k,\ell}
    = S_{k-1,\ell} + \varrho S_{k,\ell-1}
      - \varrho S_{k-1,\ell-1} +\eta_{k,\ell} ,\\
    S_{k,-1} = S_{0,\ell} = S_{0,-1} := 0,
   \end{cases}
   \qquad k \in \NN , \quad \ell \in \ZZ_+ ,
 \]
 with autoregression parameter $\varrho \in \RR$, where
 $(\eta_{k, \ell})_{k \in \NN, \ell \in \ZZ_+}$ is a set of independent standard normal random
 variables on a probability space $( \Omega, \cF, \PP )$, and $\NN$ denotes
 the set of positive integers.
The drift $\alpha_{k, \ell}$ is supposed to be an $\cF_k$-measurable random
 variable, where the filtration $( \cF_k )_{ k \in \ZZ_+ }$ is given by the trivial
 $\sigma$-algebra $\cF_0 := \{ \emptyset, \Omega \}$ and 
 \[
   \cF_k 
   := \sigma 
      ( \eta_{i,j} : \text {$1 \leq i \leq k $ \ and \ $j \geq 0$} ),
   \qquad k \in \NN .
 \]
 Let $P_{k, \ell}$ denote the price of a zero coupon bond at time $k \in \ZZ_+$
 with maturity $\ell \in \ZZ_+$ with $\ell \geq k$.
Assume that the relationship between the forward interest rates and the prices
 of a zero coupon bond is given by $P_{k,k}=1$, $k \in \NN$, and 
 \[
   P_{k, \ell + 1} = \exp \biggl\{ -\sum_{j=0}^{\ell - k} f_{k, j} \biggr\}  ,
   \qquad \text{$k, \ell \in \ZZ_+$ with $k \leq \ell$,}
 \]
 so that $P_{k, \ell + 1} = \ee^{-f_{k, \, \ell - k}} P_{k, \ell}$.
Next we consider for given positive integer $J$ a stochastic 
 discount factor process $( M_k )_{ k \in \ZZ_+ }$
 given by $M_0 := 1$ and
 \[
   M_{k+1}
   := M_k\exp{ \{- r_k \} }
      \frac{ \exp \left\{ \sum_{j=0}^{J} b_j \Delta_1 S_{k,j} \right\} 
            }{ \EE \left( \exp \left\{ \sum_{j=0}^{J} b_j \Delta_1 S_{k,j} \right\} 
                        \big| \cF_k \right)} ,
   \qquad k \in \ZZ_+ ,
 \] 
 where $r_k := f_{k,0}$ are the spot interest rate (corresponding to  
 time $k$) and
 $\mathbf{b} = (b_0, b_1, \ldots, b_J) \in \RR^{J+1}$ is the vector of the  
  market price of risk parameters. 
They play an important role in the market when determining the market prices of
 assets. 
This role is discussed in detail in \cite{GPZ}, where also the reasoning for the
 choice of the special form of the stochastic discount factors has been given. 
Note that the collection of unknown parameters we have to deal with 
  are these risk parameters, the volatility $\beta$ and the autoregression parameter  $\rho$.

We are interested only in models where arbitrage opportunities are excluded
 in the market. 
No-arbitrage property follows from a martingale condition, which is satisfied
 if the $M_k$-discounted bond price processes
 $( M_k P_{k, \ell} )_{0 \leq k \leq \ell}$ form martingales for all $\ell \in \NN$.
Using the equations resulting from the martingale condition, the drifts
 $\alpha_{k, \ell}$ disappear and we obtain
 \begin{align} \label{model}
   \begin{cases}
    f_{k,\ell} - f_{k-1,\ell+1} - \varrho ( f_{k,\ell-1} - f_{k-1,\ell} )
    = \beta \eta_{k,\ell}
      + \frac{\beta^2}{2} \sum\limits_{j=0}^{2\ell} \varrho^j
      - \beta \sum\limits_{j = \ell}^J b_j \varrho^{j - \ell} , \\
    f_{k,0} - f_{k-1,1} 
    = \beta \eta_{k,0} + \frac{\beta^2}{2}
      - \beta \sum\limits_{j = 0}^J b_j \varrho^j ,
   \end{cases}
 \end{align}
 for $k, \ell \in \NN$.
The details of the derivation of these no-arbitrage equations together with the
 role of the market discount factors can be found in \cite{GPZ}. 

The main goal of this paper is to prove 
 strong consistency and asymptotic
 normality of the joint MLE of the parameters
 $(\beta, \varrho, b_0, \dots, b_J)$ based on samples
 $(f_{k,\ell})_{1 \leq k \leq K_n, \, 0 \leq \ell \leq L_n}$, where $K_n = K n + o(n)$ 
 and
 $L_n = L n + o(n)$ as $n \to \infty$ with some $K > 0$ and $L > 0$. 
Of course, the main difficulty is 
 that the samples consist of non-independent, non-identically distributed random
 variables and moreover, no explicit formula is available for the MLE of
 $(\beta, \varrho, b_0, \dots, b_J)$.
 
 It will turn out that compared to the other two parameters 
  $\beta,$ and $\varrho,$ the 
  market price of
  risk parameters have a different asymptotic behaviour.

 When dealing with certain problems and in particular with 
  pricing derivatives, for the sake
  of convenience, many authors 
  started modelling 
  interest rate and bond markets under an 
  equivalent martingale measure. 
However,
  statistical properties of the parameter estimations usually cannot be
  discussed in that way, so that we had to work under 
  the real (objective) measure of the
  market. 
 We would like to mention that in our opinion 
  statistical tools have to be applied in
  finance for instance for pricing derivatives, 
  since in many situations the market will not be
  complete, so that one cannot work under an equivalent 
  martingale measure and one has to fit
  real date to the model. 
 Unfortunately, in the above sense 
  relatively few papers are written in finance 
  with a real statistical orientation.

Concerning the present literature we mention the following related 
 results. 
In the type of interest rate 
 framework we have investigated, there are some
 results already available for the MLE of a single parameter assuming that the
 true values of the other parameters are known.
In \cite{GPZ3}  the MLE of the volatility $\beta$ has been investigated, and
 asymptotic  normality has been obtained both in stable and in nearly unstable cases.
(A model is called stable, unstable, or explosive, if $|\varrho|<1$,
 $|\varrho|=1$, or $|\varrho|>1$, respectively. 
 In the nearly unstable case given a sequence of models with 
  corresponding autoregression parameter $\varrho_n$ we have 
  $\lim_{n \to \infty} \varrho_n = 1$.)
Volatility estimation has also been  
 studied by Peeters \cite{WP} in case of a more complicated 
 volatility structure. 
F\"ul\"op and Pap \cite{FP} tested the 
 autoregression parameter $\varrho$ both in stable and
 unstable cases, and they succeeded in proving local asymptotic normality of
 the sequence of the related statistical experiments 
 in the sense of \cite{L}. 
In a further work, in F\"ul\"op and Pap \cite{FP2}, 
 they also gave results on strong consistency of the MLE 
 estimator of the autoregressive parameter. 

The paper is organized as follows. 
In Section \ref{results} we will formulate 
 our results on consistency (Theorem \ref{stabil}) 
 and on the asymptotic normality of the joint ML parameter estimators  
 (Theorem \ref{normalitythm}).
In Section \ref{kovik} we discuss our results together with their consequences, 
 as well as some related problems and future work. 
In Appendix A we give first the derivation of the likelihood function 
 which is followed by the rigorous mathematical proofs of our main results. 
In Appendix B we collected some useful general (not model 
 specific) lemmas we apply in the 
 proofs of the main theorems.

\section{MLE and results} \label{results} 
\setcounter{theorem}{0}
\setcounter{lemma}{0}
\setcounter{proposition}{0}
\setcounter{definition}{0}
\setcounter{remark}{0}
\setcounter{corollary}{0}

In this section we present the main results on the joint 
 maximum likelihood estimators of the parameters 
 $(\beta,\varrho,b_0, \ldots , b_J)$ of the 
 model. 

Consider a sample $( f_{k,\ell} )_{ 1 \leq k \leq K, \, 0 \leq \ell \leq L }$ 
 taken from a model \eqref{model}.  
One needs first to obtain the log-likelihood function 
 which can be derived   
 based on the no-arbitrage conditions given in \cite{GPZ}. 
It has the form 
 \begin{equation} \label{likefv}
  \begin{aligned}
   \cL_{K,L}&(x_{k,\ell} 
            : 1 \leq k \leq K, \, 0 \leq \ell \leq L ; \, \beta,
            \varrho, \mathbf{b} ) 
   =-\frac{K(L+1)}{2}\log(2\pi\beta^2) \\ 
    &\quad -\frac{1}{2}\log(K!)  
      -\frac{1}{2\beta^2}
      \sum_{k=1}^K \sum_{\ell=0}^{L-1}
       \bigg(y_{k,\ell}(\varrho) 
             -\frac{\beta^2}{2}\sum_{i=0}^{2\ell}\varrho^i
             +\beta \sum_{j=\ell}^{J} b_j
                                           \varrho^{j-\ell} 
       \bigg)^2 \\[2mm]
   &\quad
     -\frac{1}{2\beta^2}
     \sum_{k=1}^{K}\frac{1}{k}
      \bigg(\ty_{k,L}(\varrho) 
            -\frac{\beta^2}{2}
             \sum_{j=1}^k \sum_{i=0}^{2(k+L-j)} \varrho^i
            +\beta \sum_{j=0}^J b_j q_{j,k,\ell} 
      \bigg)^2 , 
  \end{aligned}
 \end{equation} 
 where 
 \begin{equation}  \label{ipszilonok} 
 \begin{aligned}  
  y_{k,\ell}(\varrho)
  & := \begin{cases}  
        x_{k,\ell}-x_{k-1,\ell+1}-\varrho(x_{k,\ell-1}-x_{k-1,\ell})
        & \text{for \  $1\leq\ell\leq L-1$,}\\
        x_{k,0}-x_{k-1,1} & \text{for \  $\ell=0$,}
       \end{cases}\\
  \ty_{k,L}(\varrho)
  & := x_{k,L} - x_{0,k+L} -\varrho(x_{k,L-1} -x_{0,k+L-1})
 \end{aligned} 
 \end{equation} 
 for all $k,L \geq 1$, and $x_{0,\ell} := f_{0,\ell}$ for $\ell\geq1$. 
The derivation of the log-likelihood function is given in the 
 Appendix A in Remark \ref{likeder}. 

Unfortunately this log-likelihood function has a complicated form. 
Hence one cannot
 hope to get an explicit solutions for the estimators of all the parameters. 
We mention here that knowing the true values of some parameters, it is
 possible to give an explicit formula for the estimator(s) of the remaining
 parameter(s). 
Such a case is considered in \cite{GPZ3}, where the volatility 
 estimator is studied in a similar model. 
In general one has to use numerical procedures to maximise
 \eqref{likefv} in order to obtain the ML estimators. 
Although we do not have explicit form for the estimators, 
 the following theorems assure us that they 
 have good statistical properties (like in classical cases): 
 the first theorem is on the consistency, the second is on 
 the asymptotic normality of the joint estimators.

\begin{theorem}  \label{stabil} 
Let $H \subset \RR^{J+3}$ be a compact set such that for all
 $(\beta, \varrho, \mathbf{b}) \in H$ we have $\beta \neq 0$ and
 $\varrho \in (-1,1)$. 
Let $(\beta_0$, $\varrho_0$, $\mathbf{b}_0) \in H$ denote the true parameters,
 where we write $\mathbf{b}_0 = (b_{0,0}$, $b_{0,1}$,..., $b_{0,J})$.  
Let $K_n, L_n$, $n \in \NN$,  be positive integers such that $K_n = nK + o(n)$
 and $L_n = n L + o(n)$  as $n \to \infty$ with some $K > 0$ and $L > 0$. 
For each $n \in \NN$ let $( \hbeta_n, \hvarrho_n, \hbb_n )$ denote a maximum
 likelihood estimator of $(\beta_0, \varrho_0, \mathbf{b}_0)$ maximising the
 (log-)likelihood function over $H$. 

Then the sequence $( \hbeta_n, \hvarrho_n, \hbb_n )_{n \in \NN}$ is a strongly
 consistent estimator of $(\beta_0, \varrho_0, \mathbf{b}_0)$, i.e.,
 \begin{equation}  \label{strongly}  
  (\hbeta_n, \hvarrho_n, \hbb_n) \to (\beta_0, \varrho_0, \mathbf{b}_0) 
  \qquad\text{a.s.\ as $n \to \infty$.} 
 \end{equation}   
\end{theorem}

\begin{theorem}  \label{normalitythm} 
Under the assumptions of Theorem \ref{stabil} we have  
 \begin{equation}  \label{normi} 
  \begin{bmatrix}
   n ( \hbeta_n - \beta_0 ) \\
   n ( \hvarrho_n -\varrho_0 ) \\
   \sqrt{n} ( \hbb_n - \mathbf{b}_0 ) 
  \end{bmatrix}
  \distr \cN (0, \Lambda ),
  \qquad \text{as\  $n \to \infty$,} 
 \end{equation} 
 such that $\Lambda$ is of the form  
 \[
   \Lambda := \begin{bmatrix}  
                \Lambda_1 & \mathbf{0} \\ 
                \mathbf{0} & \Lambda_2  
              \end{bmatrix} , 
 \] 
 where 
 \[ 
   \Lambda_1 := \begin{bmatrix} 
                \sigma_{1,1} & \sigma_{1,2} \\ 
                \sigma_{2,1} & \sigma_{2,2} 
               \end{bmatrix}^{-1} 
             = \left( \sigma_{1,1} \sigma_{2,2} - \sigma_{1,2}^{2}
               \right)^{-1}   
               \begin{bmatrix} 
                \sigma_{2,2} & -\sigma_{1,2} \\ 
                -\sigma_{1,2} & \sigma_{1,1} 
               \end{bmatrix} 
 \] 
 with 
 \begin{gather}  
   \sigma_{1,1} := \frac{2KL}{\beta_0^2}  
                 + \frac{K(K+2L)}{2\left(1-\varrho_0\right)^2} , 
  \qquad 
  \sigma_{2,2} := \frac{KL}{1-\varrho_0^2} + 
                 \frac{K(K+2L)\beta_0^2}{ 
                  2\left(1-\varrho_0\right)^4} , \label{szigmak1} \\  
   \sigma_{1,2} = \sigma_{2,1} := \frac{K(K+2L)\beta_0}{
                  2\left(1-\varrho_0\right)^3} 
     \label{szigmak2}, 
 \end{gather} 
 furthermore, $\Lambda_2$ of size $(J+1) \times (J+1)$ has the form 
 \[ 
   \Lambda_2 := \frac{1}{K} 
               \begin{bmatrix} 
                 1+\varrho_0^2 & -\varrho_0 & 0 & 0 & 0 & \dots & 0 \\ 
                 -\varrho_0 & 1+\varrho_0^2 & -\varrho_0 & 0 & 0 & 
                       \dots & 0 \\ 
                 0 & -\varrho_0 & 1+\varrho_0^2 & -\varrho_0 & 0 & 
                       \dots & 0 \\ 
                 \vdots & \vdots & & & & & \vdots \\ 
                 0 & 0 & \dots & 0 & -\varrho_0 & 
                    1+\varrho_0^2 & -\varrho_0  \\ 
                 0 & 0 & \dots & \dots & 0 & -\varrho_0 & 1 
               \end{bmatrix}. 
 \] 
\end{theorem}

\section{Discussion of the results} 
\label{kovik}
\setcounter{theorem}{0}
\setcounter{lemma}{0}
\setcounter{proposition}{0}
\setcounter{definition}{0}
\setcounter{remark}{0}
\setcounter{corollary}{0}

In this paper we presented statistical results for discrete time HJM type
 forward rate models which are driven by autoregressive (AR) random fields. 
We considered some natural questions that arise when fitting such a model. 
Our aim was to examine the joint behaviour of the maximum likelihood estimator
 of all parameters of our model. 
That is, we considered the joint estimation of the AR field parameter $\varrho$, the 
 volatility parameter $\beta$ and the market price of risk parameters
 $b_0, b_1, \dots, b_J$. 

The challenge we faced was to derive good properties of the estimators
 (consistency, asymptotic normality) in a model where the observations are
 neither independent nor identically distributed. 
Furthermore, as a consequence of the complexity of the likelihood function
 there is no hope for deriving explicit solutions of the maximum likelihood 
 estimators, which complicated the task. 

Therefore, given a real market data set of forward rates, in order to fit the
 model one needs first to use numerical procedures to reach the maxima of the
 likelihood function (\ref{likefv}). 
Note that due to Theorem \ref{stabil} one can reach the maxima by 
 the use of first
 order conditions.  
On the other hand, due to the same theorem we are assured that the estimators
 are consistent. 
We also showed that joint asymptotic normality of the estimator holds like in
 the well-known cases of MLE for i.i.d.~samples (under certain conditions). 
We emphasise here that the estimators had different normalising
 factors in Theorem \ref{normalitythm}, which might be interesting for the 
 reader. 
Namely, in the normalizing factor, the market price of risk parameters differ 
 from the 'classical' (`square-root') factor 
 (of the well-known i.i.d.~cases) as the
 sample size goes to infinity. In that sense it is 
 not classical because it is not
 proportional to the reciprocal value of the square root of the sample size.
(For this notice that the sample size 
 we took in our theorem was of order $n^2$). 
Another interesting property of these risk 
 parameters is that their estimators are
 asymptotically uncorrelated from the estimators of $\beta$ and $\varrho$. 
To see this we refer to the structure 
 of the sample's Fisher information $\Sigma$ in
 Theorem \ref{anl3}.  

G\'all, Pap and Peeters 
 \cite{GPP} discussed more on the numerical problems and gave numerical results 
 of the estimations at issue. 
It was shown by the tests that even in case of small sample sizes the behaviour
 of the estimators was still very good, the estimators converged fairly fast. 
Due to this one can have good hope to fit the model well to real data. 

As we mentioned before, F\"ul\"op and Pap \cite{FP} 
 considered the separate estimation of 
 the autoregression parameter $\varrho$ both in stable and unstable cases. 
In the stable case the scaling 
 factor is $n^{-1}$, like in our case, of course.     
However, in the unstable cases the
 scaling factors are different, namely, $n^{-2}$ and $n^{-3}$. 
These scaling factors are in accordance with the Fisher information quantity
 contained in the sample. 
Based on Example 9.12 of \cite{vdV}, we expect in the explosive case 
 the sequence of the related statistical experiments to be 
 locally asymptotic mixed normal.
Finally, we note that 
 F\"ul\"op \cite{FE} gave some early numerical results on the estimation 
 of $\varrho$ as well in the above mentioned cases.

\vskip5mm
\noindent {\em Related models and future work \nopagebreak}
\vskip2mm

In G\'all, Pap and Zuijlen \cite{GPZ} a general setup has 
 been proposed for discrete time forward rate curves
 driven by random fields. 
In this paper we studied an important special case. 
However, we mention that this is certainly not the only interesting
 specification of the model one can study. 
More complicated volatility structures, other forms of market price of risk as
 well as different random fields can also be the subject of further research. 
We believe that in order to derive similar statistical results for several
 modifications of the recent model, the methods we used for the proofs will also
 work. 
For this Appendix B contains some useful tools.  
We appreciate very much the works \cite{Y1} and \cite{Y2} of Ying.
Though we did not apply directly any of Ying's specific results, his methods
 and ideas were especially fruitful during the development of the proofs of our
 main theorems.  
We note that also other methods might have been 
 also applied in order to derive the
 asymptotic results. 
Here we mention among others the excellent papers of Heijmans and Magnus 
 \cite{hm} and \cite{hm2}. 
However, the approach we have we chosen (motivated by Ying's approach) has turned out to be
 fairly appropriate and effective for our purposes. 

The asymptotic results we have found can form the basis of hypothesis tests
 that we intend to develop in our forthcoming studies. 
In this way one can hope to be able to test the goodness of fit of the model and possibly
 to compare the fits of different models. 
For model selection, information criteria might also be used. 
In our present research we are focusing on such problems. 
In that sense this paper is just the first, however the fundamental step 
 for our purposes. 
We find these problems important since, unlike in many fields of econometrics, 
 the goodnesses of fit of recently applied financial market models are often
 not justified by empirical means at all (tests, information criteria). 
They are often `just parametrised to be rich enough' so that the model produces
 (derivative) asset prices being `close enough' to market data. 
However, overparametrised models or misspecified models may occur in this way. 

Forward rate models are, of course, not only used for pricing interest rate
 derivatives. 
We hope that by finding the appropriate models and testing tools, risk
 management of firms entering to markets of bonds and interest rate related
 assets can also be better supported. 
(Here we also refer to the fact that for derivative pricing one needs not
 necessarily take our way of parameter estimation ---under the objective
 measure---, but one can alternatively use well-known calibration techniques to
 fit the models.) 
However, for many problems (e.g. risk management, goodness of fit) we suggest
 to take our approach to fit and test the model.

\section*{Acknowledgements}

The research of G. Pap was realized in the frames of
 T\'AMOP 4.2.4.\ A/2-11-1-2012-0001 ,,National Excellence Program --
 Elaborating and operating an inland student and researcher personal support
 system''.
The project was subsidized by the European Union and co-financed by the
 European Social Fund.

We are thankful to Ronald Kortram for the fruitful discussions on uniformity
 issues.

\newpage 

\section*{Appendix A: Proofs of the main results}  
\renewcommand{\theequation}{A.\arabic{equation}}
\renewcommand{\thesection}{A}
\renewcommand{\thetheorem}{A.\arabic{theorem}}
\renewcommand{\thelemma}{A.\arabic{lemma}}
\renewcommand{\theproposition}{A.\arabic{proposition}}
\renewcommand{\thedefinition}{A.\arabic{definition}}
\renewcommand{\thecorollary}{A.\arabic{corollary}}
\renewcommand{\theremark}{A.\arabic{remark}}

\setcounter{equation}{0}
\setcounter{theorem}{0}
\setcounter{lemma}{0}
\setcounter{proposition}{0}
\setcounter{definition}{0}
\setcounter{remark}{0}
\setcounter{corollary}{0}

\begin{remark}[Derivation of the likelihood function] \label{likeder} 
\textnormal{It can be seen from the main results of 
 G\'all, Pap and Zuijlen \cite{GPZ} that under the
 assumption that the common distribution of the $\eta_{i,j}$'s, for
 $i,j \in \ZZ_+$, is standard normal, the no--arbitrage criterion is equivalent
 with }
 \begin{equation}   \label{kl}
  f_{k,\ell} - f_{k-1,\ell+1} - \varrho ( f_{k,\ell-1} - f_{k-1,\ell} )
  = \beta \eta_{k,\ell}
    + \frac{\beta^2}{2} \sum_{i=0}^{2\ell} \varrho^i
    - \beta \sum_{j=\ell}^{J} b_j \varrho^{j-\ell} ,
 \end{equation}
 \textnormal{and hence }  
 \begin{equation}  \label{kl2}
  f_{k,\ell-1} - f_{k-1,\ell} 
  = \beta \sum_{i=0}^{\ell-1} \varrho^{\ell-i-1} \eta_{k,i} 
    + \frac{\beta^2}{2} \left( \sum_{i=0}^{\ell-1} \varrho^i \right)^2
    - \beta \sum_{j=0}^{J} b_j \sum_{i=0}^{j \wedge (\ell - 1)} 
       \varrho^{\ell+j-1-2i}  
 \end{equation} 
 \textnormal{for $k\geq1$, $\ell\geq1$.  
Furthermore, we have } 
 \begin{equation}   \label{kL0} 
  \begin{aligned}  
  f_{k,\ell} - f_{0,k+\ell}  
  = \sum_{n=0}^k \Bigg[ 
     &\frac{\beta^2}{2} \left( \sum_{i=0}^{k+\ell-n} \varrho^i 
                       \right)^2 
     + \beta \sum_{i=0}^{k+\ell-n} \varrho^{k+\ell-n-i} \eta_{n,i} \\ 
     &- \beta \sum_{j=0}^{J} b_j \sum_{i=0}^{j \wedge (k+\ell-n)} 
          \varrho^{k+\ell-n+j-2i} 
                 \Bigg] 
  \end{aligned} 
 \end{equation} 
 \textnormal{and } 
 \begin{equation}   \label{kL}  
  \begin{aligned}
   f_{k,\ell} - f_{0,k+\ell}
   & - \varrho ( f_{k,\ell-1} - f_{0,k+\ell-1} ) \\[2mm]
   & = \beta \sum_{j=1}^k \eta_{j,k+\ell-j}
       + \frac{\beta^2}{2} \sum_{j=1}^k \sum_{i=0}^{2(k+\ell-j)}\varrho^i
       - \beta \sum_{j=0}^J b_j q_{j,k,\ell}  
  \end{aligned}  
 \end{equation} 
 \textnormal{for $k \geq 1$, $\ell \geq 1$, 
 where }
 \begin{equation}  \label{kuk} 
     q_{j,k,\ell} := \begin{cases} 
             \sum_{n=0 \vee (j-k-\ell+1)}^{j-\ell} \varrho^{n} 
                    & \textnormal{for } j \geq \ell \\  
             0      & \textnormal{otherwise.} 
                  \end{cases} 
 \end{equation} 

\textnormal{Consider now  
 a sample $( f_{k,\ell} )_{ 1 \leq k \leq K, \, 0 \leq \ell \leq L }$
By the help of equations \eqref{kl} and \eqref{kL} one can obtain the joint
 density function of $( f_{k,\ell} )_{ 1 \leq k \leq K, \, 0 \leq \ell \leq L}$ and 
 hence the likelihood function takes the form }
 \begin{align*}
  &\bL_{K,L} (x_{k,\ell} 
              : 1\leq k\leq K, \, 0 \leq \ell \leq L ; \, \beta,
              \varrho, 
              \mathbf{b} ) 
  =\frac{1}{(2\pi\beta^2)^{(L+1)K/2}(K!)^{1/2}} \\
  &\phantom{\qquad\qquad}
   \times\exp\Bigg\{-\frac{1}{2\beta^2}
                     \sum_{k=1}^K \sum_{\ell=0}^{L-1}
                       \bigg(y_{k,\ell}(\varrho) 
                             -\frac{\beta^2}{2} \sum_{i=0}^{2\ell} \varrho^i
                             +\beta \sum_{j=\ell}^{J} b_j
                                           \varrho^{j-\ell} 
                       \bigg)^2  
   \\[2mm]
  &\phantom{\qquad\qquad\,\,\times\exp\Bigg\{}
   -\frac{1}{2\beta^2}
    \sum_{k=1}^K\frac{1}{k}
     \bigg(\ty_{k,L}(\varrho) 
           -\frac{\beta^2}{2}
            \sum_{j=1}^k \sum_{i=0}^{2(k+L-j)} \varrho^i
           + \beta \sum_{j=0}^J b_j q_{j,k,\ell}  
     \bigg)^2
    \Bigg\},
 \end{align*}
 \textnormal{where $y_{k,\ell}(\varrho)$ 
 and $\ty_{k,L}(\varrho)$ are given in \eqref{ipszilonok}. } 

\textnormal{Thus the log-likelihood function has, indeed, the form given 
 in \eqref{likefv}. 
Note that $q_{i,k,\ell}$ (see \eqref{kuk}) is bounded over a compact subset of
 $(-1,1)$ (since clearly $|q_{i,k,\ell}| \leq 1/(1-|\varrho|)$). 
Moreover recall that $q_{i,k,\ell}$ vanishes for large $\ell$. 
These facts will simplify many problems in the proofs of the results 
 on the asymptotics of the likelihood estimators. }
\end{remark}

\noindent {\bf Notation.} 
\textnormal{For simplicity, in what follows we will write } 
 \[
   \cL_n(\beta,\varrho,\mathbf{b})
   = \cL_{K_n,L_n}(f_{k,\ell}: 1  \leq k \leq K_n, 0 \leq \ell \leq L_n ;  
                          \beta, \varrho, \mathbf{b} )  
 \] 
 \textnormal{and } 
 \begin{equation*}  
 \begin{aligned} 
    \partial_1^{i_1} \partial_2^{i_2} \partial_{j_1}^{i_3}  
    \partial_{j_2}^{i_4}  
   &\cL_n(\beta,\varrho ,\mathbf{b}) \\ 
   &= \frac{\partial^{i_1} \partial^{i_2} \partial^{i_3}
    \partial^{i_4}  
     \cL_{K_n,L_n}(
     x_{k,\ell}:1\  \leq k \leq K_n, 0 \leq \ell \leq L_n; 
     \beta, \varrho , \mathbf{b})}{ 
         \partial \beta^{i_1} \partial \varrho^{i_2} 
         \partial b_{j_1-3}^{i_3} 
         \partial b_{j_2-3}^{i_4}  
                             } 
    \bigg|_{x_{k,\ell}=f_{k,\ell}} \; ,  
 \end{aligned} 
 \end{equation*}  
 \textnormal{where $i_1, i_2, i_3, i_4, j_1, j_2$ are non-negative integers and
 $3 \leq j_i \leq J+3$ for $i = 1, 2$. 
Furthermore, } 
 \[
   \pmb{\partial_3} \cL_n(\beta, \varrho,\mathbf{b})
   = \begin{bmatrix} 
      \partial_3 \cL_n(\beta, \varrho,\mathbf{b}) \\ 
      \partial_4 \cL_n(\beta, \varrho, \mathbf{b}) \\  
      \vdots \\ 
      \partial_{J+3} \cL_n(\beta, \varrho, \mathbf{b}) 
     \end{bmatrix} . 
 \]


\noindent 
\textbf{Proof of Theorem \ref{stabil}.} 
First we show strong consistency of $\hbeta_n,$ and $\hvarrho_n$. 
For this, the aim of the following discussion is to derive an asymptotic
 expansion for the sequence of random variables 
 \begin{equation*}   
  \cL_n ( \beta, \varrho, \mathbf{b} )
   = \cL_{K_n,L_n}(f_{k,\ell}
                   : 1 \leq k \leq K_n, \, 0 \leq \ell \leq L_n ; \,
                     \beta, \varrho, \mathbf{b}) ,
  \qquad n \geq 1. 
 \end{equation*} 
We have
 \begin{equation}  \label{hn} 
  \begin{aligned}
  \cL_n ( \beta, \varrho, \mathbf{b} )
  &= -\frac{K_n(L_n+1)}{2}\log(2\pi\beta^2)-\frac{1}{2}\log(K_n!) \\ 
  &\phantom{=\:}   
     -\frac{1}{2\beta^2}
      \sum_{k=1}^{K_n} \sum_{\ell=0}^{L_n-1} 
          \xi_{k,\ell}^2 ( \beta, \varrho, \mathbf{b} ) \\
  &\phantom{=\:}
     -\frac{1}{2\beta^2}
      \sum_{k=1}^{K_n} k^{-1}
      \left(\sum_{j=1}^k 
         \xi_{j,k+L_n-j} ( \beta, \varrho, \mathbf{b} ) 
      \right)^2, 
  \end{aligned} 
 \end{equation} 
 where
 \begin{equation}  \label{xidef} 
  \xi_{k,\ell} (\beta, \varrho, \mathbf{b}) 
  := g_{k,\ell}(\varrho) 
     - \frac{\beta^2}{2} \sum_{i=0}^{2\ell} \varrho^i
     + \beta \sum_{j=\ell}^{J} b_{j} \varrho^{j-\ell} 
 \end{equation} 
 with 
 \begin{align*}
  g_{k,\ell}(\varrho) 
  & := \begin{cases}
        f_{k,\ell}-f_{k-1,\ell+1}-\varrho(f_{k,\ell-1}-f_{k-1,\ell})
        & \text{for \  $\ell \geq 1$,}\\
        f_{k,0}-f_{k-1,1} & \text{for \  $\ell=0$,}
       \end{cases}
 \end{align*}
 for all \  $k \geq 1$. 
Since we have 
 $g_{k,\ell}(\varrho)
  = g_{k,\ell}(\varrho_0) + ( \varrho_0 - \varrho ) ( f_{k,\ell-1} - f_{k-1,\ell} )$, 
 $\ell \geq 1$,  
 by applying formula \eqref{kl2} we obtain
 \begin{equation}  \label{kszi} 
  \begin{aligned}
  \xi_{k,\ell} ( \beta, \varrho, \mathbf{b} ) 
  &= \beta_0\eta_{k,\ell}
     + \frac{\beta_0^2}{2} \sum_{i=0}^{2\ell} \varrho_0^i
     - \frac{\beta^2}{2} \sum_{i=0}^{2\ell} \varrho^i
     - \beta_0 \sum_{j=\ell}^{J} b_{0,j} \varrho_0^{j-\ell} 
     + \beta \sum_{j=\ell}^{J} b_{j} \varrho^{j-\ell}  \\
  &\phantom{\quad}
     + ( \varrho_0 - \varrho )
       \Bigg[ \beta_0 \sum_{i=0}^{\ell-1} \varrho_0^{\ell-i-1} \eta_{k,i} 
              + \frac{\beta_0^2}{2}
                \left( \sum_{i=0}^{\ell-1} \varrho_0^i \right)^2 \\ 
  &\phantom{\quad + ( \varrho_0 - \varrho ) \Bigg[ \;} 
              -\beta_0 \sum_{j=0}^{J} b_{0,j} 
                     \sum_{i=0}^{j \wedge \ell} \varrho_0^{\ell+j-2i}         
        \Bigg] .
  \end{aligned}
 \end{equation} 
We have
 \begin{align*}
  \EE \xi_{k,\ell} ( \beta, \varrho, \mathbf{b} )
  &= \frac{\beta_0^2}{2} \sum_{i=0}^{2\ell} \varrho_0^i
     - \frac{\beta^2}{2} \sum_{i=0}^{2\ell} \varrho^i
     - \beta_0 \sum_{j=\ell}^{J} b_{0,j} \varrho_0^{j-\ell} 
     + \beta \sum_{j=\ell}^{J} b_{j} \varrho^{j-\ell}  \\
  &\phantom{\quad}
     + ( \varrho_0 - \varrho )
       \left[ \frac{\beta_0^2}{2}
              \left( \sum_{i=0}^{\ell-1} \varrho_0^i \right)^2
              -\beta_0 \sum_{j=0}^{J} b_{0,j} 
                     \sum_{i=0}^{j \wedge \ell} \varrho_0^{\ell+j-2i}    
       \right] 
   \to m( \beta, \varrho )
 \end{align*}
 as $\ell \to \infty$, where
 \begin{equation}  \label{emi} 
  m ( \beta, \varrho ):= \frac{\beta_0^2}{2(1-\varrho_0)} 
       - \frac{\beta^2}{2(1-\varrho)} 
       + \frac{(\varrho_0-\varrho)\beta_0^2}{2(1-\varrho_0)^2} .
 \end{equation} 
Hence $\sup_{k,\ell} | \EE \xi_{k,l} ( \beta, \varrho, \mathbf{b} )| < \infty$.
Moreover,
 \[
   \Var \xi_{k,\ell} ( \beta, \varrho, \mathbf{b} )
   = \beta_0^2 
     \left( 1 + ( \varrho_0 -\varrho )^2
                \sum_{i=0}^{\ell-1} \varrho_0^{2(\ell-i-1)} \right)
   \to \sigma^2 ( \varrho )
 \]
 as $\ell \to \infty$, where
 \begin{equation}  \label{szigi}  
  \sigma^2( \varrho )
  := \beta_0^2 \left( 1 + \frac{(\varrho_0-\varrho)^2}{1-\varrho_0^2} \right) .
 \end{equation} 
Hence $\sup_{k,\ell} | \Var \xi_{k,l} ( \beta, \varrho, \mathbf{b} )| < \infty$.
Since $\xi_{k,\ell} ( \beta, \varrho, \mathbf{b} )$ has a normal distribution for
 all $k \geq 1$, $\ell \geq 0$, we conclude 
 $\sup_{k,\ell} \EE \xi_{k,l}^8 ( \beta, \varrho, \mathbf{b} ) < \infty$. 
Furthermore, 
 \[
   n^{-2} \sum_{k=1}^{K_n} \sum_{\ell=0}^{L_n-1} \EE 
     \xi_{k,\ell}^2 ( \beta, \varrho, \mathbf{b} ) 
   \to KL ( \sigma^2( \varrho ) 
           + m^2( \beta, \varrho ) )
 \]
 as $n \to \infty$. 
Obviously the sets
 $\{ \xi_{k,\ell}( \beta, \varrho, \mathbf{b} ) : \ell \in \NN \}$,
 $k \in \NN$, are independent, hence by Lemma \ref{lem1} we obtain
 \begin{equation}  \label{egyes} 
  n^{-2} \sum_{k=1}^{K_n} \sum_{\ell=0}^{L_n-1} 
   \xi_{k,\ell}^2 ( \beta, \varrho, \mathbf{b} ) 
  \to KL ( \sigma^2( \varrho ) + m^2( \beta, \varrho ) )
  \qquad \text{a.s. as \  $n \to \infty$.}
 \end{equation}
Clearly $\{ \xi_{j,k+L_n-j} ( \beta, \varrho, \mathbf{b} ) : 1 \leq j \leq k \}$
 are independent for all $k, n \in \NN$, hence
 \[
   \EE \left(\sum_{j=1}^k \xi_{j,k+L_n-j} ( \beta, \varrho, \mathbf{b} )
       \right)^2
   = \sum_{j=1}^k \Var\xi_{j,k+L_n-j}( \beta, \varrho, \mathbf{b} )
     + \left(\sum_{j=1}^k \EE \xi_{j,k+L_n-j} ( \beta, \varrho, \mathbf{b} )
       \right)^2 .
 \] 
Applying the above formulas for $\EE \xi_{k,\ell} ( \beta, \varrho, \mathbf{b} )$
 and $\Var \xi_{k,\ell} ( \beta, \varrho, \mathbf{b} )$ it follows that
 \[
   n^{-2} \sum_{k=1}^{K_n} k^{-1} 
   \EE \left(\sum_{j=1}^k \xi_{j,k+L_n-j} ( \beta, \varrho, \mathbf{b} )
   \right)^2
   \to \frac{K^2m^2( \beta, \varrho )}{2} 
 \]
 as $n \to \infty$, hence by Lemma \ref{lem2} we obtain
 \begin{equation}  \label{kettes}
  n^{-2} \sum_{k=1}^{K_n} k^{-1} \left(\sum_{j=1}^k 
    \xi_{j,k+L_n-j}( \beta, \varrho, \mathbf{b} ) \right)^2
  \to \frac{K^2m^2( \beta, \varrho )}{2} 
  \qquad \text{a.s. as $n \to \infty$.}
 \end{equation}
Now, equations \eqref{egyes} and \eqref{kettes} lead us to 
 \begin{align}  
   & \cL_n ( \beta_0, \varrho_0, \mathbf{b}_0 )
    - \cL_n ( \beta, \varrho, \mathbf{b} )  
     \nonumber \\[2mm]
   & = \frac{KLn^2}{2}
        \left( \frac{\beta_0^2}{\beta^2} - 1
               - \log \frac{\beta_0^2}{\beta^2} \right) 
        + \frac{KL\left(\varrho_0-\varrho\right)^2\beta_0^2n^2}
               {2\beta^2\left(1-\varrho_0^2\right)} \label{fosor}  \\[2mm] 
   &\phantom{=\:}
        + \frac{K(K+2L)n^2}{16\beta^2}
          \left( \frac{\beta_0^2\left(\varrho_0-\varrho\right)}
                      {\left(1-\varrho_0\right)^2}
                 + \frac{\beta_0^2}{1-\varrho_0} 
                 - \frac{\beta^2}{1-\varrho} \right)^2
        + o (n^2) \nonumber  
   \qquad \text{a.s. as $n \to \infty$.}
  \end{align}  
Furthermore, notice that \eqref{fosor} holds uniformly in   
 $( \beta, \varrho, \mathbf{b})$ over $H$, due to the special from of the
 likelihood function. 
We show the details of the proof of uniformity in Remark \ref{uniform}. 

For a fixed $n$, one can now consider a maximum likelihood estimator of 
 $( \beta_0, \varrho_0, \mathbf{b}_0 )$, say $( \hbeta_n, \hvarrho_n, \hbb_n )$,
 which is the maximiser of $\cL_n ( \beta, \varrho, \mathbf{b} )$ over $H$. 
Hence, after replacing $( \beta, \varrho, \mathbf{b} )$ by 
 $( \hbeta_n, \hvarrho_n, \hbb_n )$ in \eqref{fosor} one can easily see that the
 left hand side is non-positive with probability one, that is
 a.s.~$\cL_n(\beta_0, \varrho_0, \mathbf{b}_0 )
       - \cL_n ( \hbeta_n, \hvarrho_n, \hbb_n ) \leq 0$. 
On the other hand the leading terms of the right hand side of \eqref{fosor} are
 non-negative and at least one of them is positive if
 $( \hbeta_n, \hvarrho_n) \neq (\beta_0, \varrho_0)$. 
Therefore, as $n \to \infty$, equation \eqref{fosor} can be kept only if
 $( \hbeta_n, \hvarrho_n) \to (\beta_0, \varrho_0)$ a.s., since the right hand
 side of \eqref{fosor} would not tend to $0$ as $n \to \infty$ for
 $\omega \in \Omega$ if we had 
 $\big( \hbeta_n(\omega), \hvarrho_n(\omega) \big)
  \not \to (\beta_0, \varrho_0)$. 
That is, the strong consistency of the maximum likelihood estimators of
 $(\beta,\varrho)$ holds. 

Now we turn to showing strong consistency of $\hbb_n$. 
Consider the system of equations determined by the first order conditions
 $\partial_{j+3} \cL_n(\hbeta_n, \hvarrho_n, \hbb_n)=0$ for $j=0, 1, \ldots , J$. 
For large $n$ (e.g. for $L_n > J$, for this recall that due the remark made on
 the vanishing $q_{j,k,\ell}$'s at the end of Remark \ref{likeder} only the
 second line of the right hand side of \eqref{likefv} will contain the market
 price of risk parameters) we can rewrite this system of equations in the
 simple form (see \eqref{hn} and \eqref{xidef}) 
 \begin{equation*} 
  \sum_{k=1}^{K_n} \sum_{i=0}^j \xi_{k,i} ( \hbeta_n, \hvarrho_n, \hbb_n ) 
    \hvarrho^{\; j-i} = 0  \quad \text{for } j=0, 1,  \ldots , J, 
 \end{equation*} 
 which can be reduced to  
 \begin{equation}  \label{bkonz}  
   \sum_{k=1}^{K_n} \xi_{k,j} ( \hbeta_n, \hvarrho_n, \hbb_n ) = 0 
    \quad \text{for } j=0, 1, \ldots , J.  
 \end{equation} 
Now, taking \eqref{bkonz} for $j=J$ we obtain 
 \begin{equation}  \label{bj} 
    \sum_{k=1}^{K_n} \left[ \beta_0 \eta_{k,J} 
    + \frac{\beta_0^2}{2} \sum_{i=0}^{2J} \varrho_0^i 
    - \frac{\hbeta_n^2}{2} \sum_{i=0}^{2J} \hvarrho_n^{\; i} 
    - \beta_0 b_{0,J} + \hbeta_n \; \hb_{n,J}
    + \left( \varrho_0 - \hvarrho_n \right) c^{(2)}_{k,J} \right] 
    = 0, 
 \end{equation}  
 where for $k, \ell \in \ZZ_+$, $k > 0$ we write 
 \begin{equation}  \label{c2kl}  
   c^{(2)}_{k,\ell} := 
    \beta_0 \sum_{i=0}^{\ell-1} \varrho_0^{\ell-i-1} \eta_{k,i} 
   + \frac{\beta_0^2}{2}
      \left( \sum_{i=0}^{\ell-1} \varrho_0^i \right)^2 
   - \beta_0 \sum_{j=0}^{J} b_{0,j} 
             \sum_{i=0}^{j \wedge \ell} \varrho_0^{\ell+j-2i}  .  
 \end{equation} 
Notice that the random variable $c^{(2)}_{k,\ell}$ does not depend on 
 $( \hbeta_n, \hvarrho_n, \hbb_n )$ and $\{c^{(2)}_{k,\ell}\}_{k \geq 0}$ are
 i.i.d.~for a fixed $\ell \in \ZZ_+$. 
Reordering \eqref{bj} we obtain 
 \[
   \hb_{n,J} - b_{0,J} = 
      \frac{\beta_0-\hbeta_n}{\hbeta_n} b_{0,J} 
      - \frac{\beta_0^2}{2\hbeta_n} \sum_{i=0}^{2J} \varrho_0^i 
      + \frac{\hbeta_n}{2} \sum_{i=0}^{2J} \hvarrho_n^{\; i} 
      - \frac{1}{K_n \hbeta_n} \sum_{k=1}^{K_n} 
         \left[ \beta_0 \eta_{k,J} 
           + \left( \varrho_0 - \hvarrho_n \right) c^{(2)}_{k,J} 
         \right] . 
 \] 
Hence, by the SLLN and the consistency of $(\hbeta_n, \hvarrho_n)$ we obtain
 that $\hb_{n,J} \to b_{0,J}$ a.s.\ as $n\to\infty$, i.e.~$\hb_{n,J}$ is strongly
 consistent. 
In a similar way, recursively we can obtain the consistency of
 $\hb_{n,J-1}, \hb_{n,J-2}, \ldots , \hb_{n,1}$. 
Indeed, given the consistency of $\hb_{n,J-1}, \hb_{n,J-2}, \ldots , \hb_{n,j+1}$,
 consider again \eqref{bkonz} from which we can obtain 
 \begin{equation*} 
  \begin{aligned} 
    \hb_{n,j} - b_{0,j} = 
      &\frac{\beta_0-\hbeta_n}{\hbeta_n} b_{0,j} 
      + \sum_{i=j+1}^{J} \left( \frac{\beta_0}{\hbeta_n} 
               \; b_{0,i} \varrho_0^{i-j} 
             - \hb_{n,i} \; \hvarrho_n^{\; i-j} 
                         \right) \\ 
      &- \frac{\beta_0^2}{2\hbeta_n} \sum_{i=0}^{2j} \varrho_0^i 
      + \frac{\hbeta_n}{2} \sum_{i=0}^{2j} \hvarrho_n^{\; i}         
      - \frac{1}{K_n \hbeta_n} \sum_{k=1}^{K_n} \left[ \beta_0 \eta_{k,j} 
      + \left( \varrho_0 - \hvarrho_n \right) c^{(2)}_{k,j} 
                                       \right] ,  
  \end{aligned}      
 \end{equation*} 
 from which the consistency of $\hb_{n,j}$ follows and thus the proof of Theorem
 \ref{stabil} is complete.
\proofend

\begin{remark}[Uniformity in (\ref{fosor})]  \label{uniform} 
\textnormal{In the derivation of \eqref{fosor} in fact we have shown that for
 any fixed point $(\beta, \varrho, \pmb{b}) \in H$ we have}
 \[
   n^{-2} \left( \cL_n(\beta, \varrho, \pmb{b}) 
      + \frac{1}{2} \log (K_n !) 
            \right) 
     \to 
     A(\beta, \varrho) ,
     \qquad 
     a.s., 
 \] 
 \textnormal{where, recalling notations \eqref{emi} and \eqref{szigi}, 
  the (deterministic) function $A$ is defined as} 
  \[
    A(\beta, \varrho) = -\frac{KL}{2} \log \left(2\pi \beta^2\right) 
      - \frac{KL \left( \sigma^2(\varrho)+m^2(\beta , \varrho) 
                 \right) }{2\beta^2}  
      - \frac{K^2 m^2(\beta , \varrho) }{4 \beta^2} , 
 \] 
 \textnormal{for $(\beta, \varrho) \in \RR^2$, $\beta \neq 0$.} 

\textnormal{Now, introducing the notations 
 $ c^{(1)}_{\ell}(\beta, \varrho, \pmb{b}) 
   := \frac{\beta^2 \varrho^{2\ell + 1}}{2(1-\varrho)} 
      + \beta \sum_{j=\ell}^{J} b_{j} \varrho^{j-\ell}$ 
 and $c^{(3)}_{\ell}
      := \frac{\beta_0^2}{2} \sum_{i=0}^{2\ell} 
         \varrho_0^i - \beta_0 \sum_{j=\ell}^{J} b_{0,j} \varrho^{j-\ell}$ 
 we can rewrite \eqref{kszi} as} 
 \begin{equation}  \label{ujalak} 
    \xi_{k,\ell} ( \beta, \varrho, \mathbf{b} )= \beta_0 \eta_{k,\ell} 
     -\frac{\beta^2}{2\left( 1-\varrho \right)}  
     + c^{(1)}_{\ell}(\beta, \varrho, \pmb{b}) 
     + ( \varrho_0 - \varrho ) c^{(2)}_{k,\ell} + c^{(3)}_{\ell} , 
 \end{equation}  
 \textnormal{where $c^{(2)}_{k,\ell}$ is given in \eqref{c2kl}. 
Notice that $c^{(2)}_{k,\ell}$ is a random variable, $c^{(3)}_{\ell}$ is a constant
 and these latter two terms depend only on $(\beta_0, \varrho_0, \pmb{b}_0)$
 but not on $(\beta, \varrho, \pmb{b})$. 
In this way of writing $\xi_{k,\ell}$ we have displayed only the parts which
 depend on the parameters $(\beta, \varrho, \pmb{b})$. 
We can see that this dependence is relatively simple and, say, fairly separated
 from the random parts.}

\textnormal{Now take the square of $\xi_{k,\ell}( \beta, \varrho, \mathbf{b} )$ 
 based on \eqref{ujalak} 
 and substitute it in \eqref{hn}. 
In the followings we will consider the terms we obtain in the square of
 $\xi_{k,\ell}( \beta, \varrho, \mathbf{b} )$. 
We mention that in \eqref{ujalak} we displayed the term 
 $-\frac{\beta^2}{2\left( 1-\varrho \right)}$ rather than embedding it in
 $c^{(1)}_{\ell}(\beta, \varrho, \pmb{b})$. 
The reason for that was that the terms (of the log-likelihood function) which
 contain $c^{(1)}_{\ell}(\beta, \varrho, \pmb{b})$ will be shown to vanish
 uniformly as $n \to \infty$ unlike the terms containing
 $-\frac{\beta^2}{2\left( 1-\varrho \right)}$.}\par 
\textnormal{By the application of Lemmas \ref{lem1}, \ref{lem1_2}, \ref{lem2} 
 and their corollaries (see Appendix B) 
 we can easily see that for $m=1,2$ the terms} 
 \begin{gather*}  
  n^{-2} \sum_{k=1}^{K_n} \sum_{\ell=0}^{L_n-1} 
     \left( \eta_{k,\ell} \right)^m, \qquad 
  n^{-2} \sum_{k=1}^{K_n} k^{-1} 
     \left( \sum_{j=1}^k \eta_{j,k+L_n-j}
     \right)^m,  \\  
  n^{-2} \sum_{k=1}^{K_n} \sum_{\ell=0}^{L_n-1} 
     \left( c^{(2)}_{k,\ell} 
     \right)^m, \qquad 
  n^{-2} \sum_{k=1}^{K_n} k^{-1} 
     \left( \sum_{j=1}^k c^{(2)}_{j,k+L_n-j}
     \right)^m,  \\ 
  n^{-2} \sum_{k=1}^{K_n} \sum_{\ell=0}^{L_n-1} 
    \eta_{k,\ell} c^{(2)}_{k,\ell},  \qquad 
  n^{-2} \sum_{k=1}^{K_n} k^{-1} 
    \eta_{j,k+L_n-j} c^{(2)}_{j,k+L_n-j}  
 \end{gather*} 
 \textnormal{all have an almost sure limit. 
Therefore, let $\Gamma_{\beta_0, \varrho_0, \pmb{b}_0} \subset \Omega$ denote the set
 over which these terms all converge to the their limits (given by Lemmas
 \ref{lem1}, \ref{lem1_2}, \ref{lem2} and their corollaries). 
Thus, $\PP(\Gamma_{\beta_0, \varrho_0, \pmb{b}_0})=1$. 
We will show that the uniformity of the almost sure convergence at issue is
 fulfilled over this set.}\par 
\textnormal{Next consider the terms we obtain in \eqref{hn} after taking the
 square of $\xi_{k,\ell}( \beta, \varrho, \mathbf{b} )$ based on \eqref{ujalak}
 which contain $c^{(1)}_{\ell}$. 
According to our assumptions $\sup_{(\beta, \varrho, \pmb{b}) \in H} | \varrho | < 1$. 
Hence observe that $\left| \sum_{\ell=0}^{L_n-1} c^{(1)}_{\ell} \right|$ and
 $\left| k^{-1} \sum_{j=1}^k c^{(1)}_{k+L_n-j} \right|$ are both bounded above for
 $k, n$. 
Therefore for $m=1,2$ the terms} 
 \[
   n^{-2} \sum_{k=1}^{K_n} \sum_{\ell=0}^{L_n-1} 
    \left( c^{(1)}_{\ell} 
     \right)^m \qquad \text{and}  
    \qquad  
   n^{-2} \sum_{k=1}^{K_n} k^{-1} 
     \left( \sum_{j=1}^k c^{(1)}_{k+L_n-j}
     \right)^m 
 \] 
 \textnormal{vanish uniformly over $H$.}\par 
\textnormal{Furthermore, the cross product terms containing $c^{(1)}_{\ell}$ all
 vanish uniformly in $H$. 
One can see this by applying the Cauchy-Schwartz inequality. 
For instance,} 
 \[
   \left| n^{-2} \sum_{k=1}^{K_n} \sum_{\ell=0}^{L_n-1} 
      c^{(1)}_{\ell} c^{(2)}_{k,\ell} \right| \leq 
       \left[ n^{-2} \sum_{k=1}^{K_n} \sum_{\ell=0}^{L_n-1} 
              \left( c^{(1)}_{\ell} \right)^2 \right]^{1/2}
       \left[ n^{-2} \sum_{k=1}^{K_n} \sum_{\ell=0}^{L_n-1} 
              \left( c^{(2)}_{k,\ell} \right)^2 \right]^{1/2}
      \to 0 
 \] 
 \textnormal{as $n \to \infty$ uniformly in $H$.}\par 
\textnormal{It is easy to check that the remaining terms we obtained in
 \eqref{hn} also converge almost surely and uniformly over $H$.}\par 
\textnormal{Summarising the above results we obtain that on the one hand} 
 \begin{equation}  \label{univegre} 
    n^{-2} \left( \cL_n(\beta, \varrho, \pmb{b}) 
         + \frac{1}{2} \log (K_n !) \right)  
     \to 
     A(\beta, \varrho) ,
     \quad \forall \omega \in \Gamma_{\beta_0, 
       \varrho_0, \pmb{b}_0}, \; \; \forall 
       (\beta, \varrho, \pmb{b}) \in H ,  
 \end{equation} 
 \textnormal{and on the other hand the uniformity of the convergences detailed 
 in the last paragraphs imply that \eqref{univegre} holds uniformly over $H$,
 which means that the almost sure expansion \eqref{fosor} holds uniformly in
 $(\beta, \varrho, \pmb{b}) \in H$, indeed.}
\end{remark}


{\bf Proof of Theorem \ref{normalitythm}.} 
We apply again Taylor's expansion for the gradient vector of 
 $\cL_n(\beta, \varrho, b )$ up to order 2. 
Write 
 \begin{equation}  \label{felbontas}  
   \begin{bmatrix}  
       n^{-1} \partial_1 \cL_n(\hbeta_n,\hvarrho_n,\hbb_n) \\  
       n^{-1} \partial_2 \cL_n(\hbeta_n,\hvarrho_n,\hbb_n) \\  
       n^{-1/2} \pmb{\partial_3} \cL_n(\hbeta_n,\hvarrho_n,\hbb_n) \\ 
   \end{bmatrix} 
  -\begin{bmatrix}  
       n^{-1} \partial_1 \cL_n(\beta_0,\varrho_0,\mathbf{b}_0) \\  
       n^{-1} \partial_2 \cL_n(\beta_0,\varrho_0,\mathbf{b}_0) \\  
       n^{-1/2} \pmb{\partial_3} \cL_n(\beta_0,\varrho_0,\mathbf{b}_0) \\ 
   \end{bmatrix}  =  
  \left( \mathbf{A}_n + \mathbf{B}_n 
  \right) 
  \begin{bmatrix} n ( \hbeta - \beta_0 )  \\  
                  n ( \hvarrho - \varrho_0 ) \\  
                  \sqrt{n} ( \hbb - \mathbf{b}_0 ) \\ 
  \end{bmatrix} 
 \end{equation} 
 where $\mathbf{A}_n$ and $\mathbf{B}_n$ are $(J+3) \times (J+3)$ matrices
 defined as follows. 
Write $\mathbf{A}_n := (a_{i,j}^n)_{i,j=1,\ldots , J+3}$ and define 
 \[
   a_{i,j}^n :=  n^{l_{i,j}} \partial_i \partial_j 
    \cL_n (\beta_0, \varrho_0, \mathbf{b}_0),
 \] 
 where 
 \[
   l_{i,j} := \begin{cases} 
                -2 & \text{if } i \vee j \leq 2 \\  
                -1 & \text{if } i \wedge j \geq 3 \\  
                -3/2 & \text{otherwise.} 
               \end{cases}  
 \] 
Denoting the $i$th row of $\mathbf{B}_n$ by $B_n^i$ we will write it in the
 form $B_n^i = D_n^{\top} R_n^i$, where the superscript ${ }^{\top}$ denotes the
 transposed, 
 \[
   D_n^{\top}
   := ( \hbeta_n - \beta_0,\; \hvarrho_n - \varrho_0,\;  
                   \hbb_n - \mathbf{b}_0),
 \] 
 \[
   R_n^i := \left( r^{n,i}_{j_1,j_2} \right)_{j_1,j_2=1,2,\ldots,J+3}
 \] 
 and 
 \[
   r^{n,i}_{j_1,j_2} := \frac{1}{2} n^{l_{i,j_2}}\partial_i \partial_{j_1} 
     \partial_{j_2} \cL_n ( \tilde{\beta}, \tilde{\varrho}, 
     \tilde{\mathbf{b}} )  
 \] 
 with appropriate $(\tilde{\beta}, \tilde{\varrho}, \tilde{\mathbf{b}} )$  
 taking values ---coordinate-wise--- between
 $(\beta_0, \varrho_0, \mathbf{b}_0 )$ and $(\hbeta_n, \hvarrho_n, \hbb_n)$.  

Under the assumptions of Theorem \ref{stabil} we will need the following lemmas
 on the terms introduced in (\ref{felbontas}).  
The proofs of these lemmas follow this proof. 
\begin{lemma}  \label{anl3} 
 Under the assumptions of Theorem \ref{stabil} we have 
 \begin{equation}  \label{ll3} 
   \begin{bmatrix} n^{-1} \partial_1 \cL_n(\beta_0,\varrho_0,
                                      \mathbf{b}_0) \\  
                   n^{-1} \partial_2 \cL_n(\beta_0,\varrho_0,
                                      \mathbf{b}_0) \\  
                   n^{-1/2} \pmb{\partial_3} \cL_n(\beta_0,\varrho_0,
                                      \mathbf{b}_0) 
   \end{bmatrix}
  \distr \cN(0,\Sigma)  , 
 \end{equation}  
 where $\Sigma = \left( \sigma_{i,j} \right)_{i,j=1, \ldots, J+3}$ with $\sigma_{i,j}$
 for $i \vee j \leq 2$ defined in \eqref{szigmak1} and \eqref{szigmak2} in
 Theorem \ref{normalitythm}, 
 \begin{equation*}    
   \sigma_{i+3,j+3} = K \sum_{k=0}^{i \wedge j} \varrho^{i+j-2k} 
        \quad \text{for $i, j = 0, 1, \ldots , J,$} 
 \end{equation*} 
 and the remaining entries of $\Sigma$ are zero. 
\end{lemma}   

\begin{lemma}  \label{anl1} 
 \[
   \mathbf{A}_n \to - \Sigma \qquad \text{a.s. as } n \to \infty.
 \] 
\end{lemma} 
\begin{lemma}  \label{anl2} 
 For $i=1,2,\ldots,J+3$ we have 
 \[
   B_n^i \stoch 0 .
 \]
\end{lemma} 

Clearly, the first term on the left hand side of \eqref{felbontas} tends to
 zero almost surely, in fact it takes value 0 a.s.~for large $n$ due to Theorem
 \ref{stabil}. 
Hence, by Slutsky's Lemma the limit distribution of the left hand side of
 \eqref{felbontas} is $\cN(0,\Sigma)$, which is given in Lemma \ref{anl3}. 
Lemma \ref{anl1} and Lemma \ref{anl2} together with Slutsky's Lemma give 
 \begin{equation*}  
    A_n + B_n \stoch -\Sigma .  
 \end{equation*}  
Having these asymptotic results and recalling \eqref{felbontas} we can apply
 Lemma \ref{anl4} to obtain \eqref{normi}. 
For this note that $\Lambda = \Sigma^{-1}$ where $\Sigma$ is given in Lemma
 \ref{anl3}. 
Thus, for the proof of Theorem \ref{normalitythm} there remains to prove Lemma
 \ref{anl3}, Lemma \ref{anl1} and Lemma \ref{anl2}. 
\proofend

\noindent {\bf Proof of Lemma \ref{anl3}.} 
First we will show that 
 \begin{equation*}  
   \mathbf{Z}_n := 
    \begin{bmatrix} n^{-1} \partial_1 \cL_n(\beta_0,\varrho_0,\pmb{b}_0) \\  
                   n^{-1} \partial_2 \cL_n(\beta_0,\varrho_0,\pmb{b}_0) \\  
                   n^{-1/2} \pmb{\partial}_3 
                      \cL_n(\beta_0,\varrho_0,\pmb{b}_0) \\  
    \end{bmatrix}
 \end{equation*} 
 can be considered as a martingale with respect to an appropriate filtration.  
Namely, rewriting the terms of $\mathbf{Z}_n$ in an appropriate order we will
 obtain the form 
 \begin{equation*} 
  \mathbf{Z}_n = \sum_{m=1}^{K_n(L_n+1)} M_{m}^{(n)} , 
 \end{equation*} 
 where $M_{m}^{(n)}$ is defined below. 
The idea of reordering the terms is simple: starting with and fixing $k=0$ we
 increase $\ell$ step by step (as $m$ increases) from $0$ to $L_n$. 
When $\ell=L_n$ is reached after $L_n+1$ steps than we consider the next value
 of $k$ ($=1$) and we take again the possible values of $\ell$ from $0$ to
 $L_n$. 
We continue this as long as $k=K_n$ is reached. 
Thus the number of summands is $K_n(L_n+1)$. 
This means that in each step a new $\eta_{k,\ell}$ will occur in the martingale
 sum (which is independent of the previous terms). 
The case $\ell = L_n$ is a little bit special, since it involves a number of
 new variables, namely, $\eta_{k,L_n}$, $\eta_{k-1,L_n+1}$, ... $\eta_{1,L_n+k-1}$ 
 (which are also independent of the previous terms). 
 
Let us turn now to the rigorous definition of the martingale difference
 $M_{m}^{(n)}$, and the corresponding filtration. 
Fix $n\in\NN$ and notice that for any positive integer $m$ there exist uniquely
 determined integers $k_m,\ell_m$ such that $m=(k_m-1)(L_n+1)+\ell_m+1$ with
 $0<k_m$ and $0\leq \ell_m \leq L_n$. 
We remark that $k_m$ and $\ell_m$ depend on $n$ as well, however, for
 simplicity we omit to denote their dependence on $n$. 
Now, define $\cG_0^{(n)} := \{ \emptyset , \Omega \}$ and 
 \begin{equation*} 
  \cG_m^{(n)} : = \begin{cases} 
    \sigma\left\{\cG_{m-1}^{(n)} \bigcup \sigma 
                     \left\{\eta_{k_m,\ell_m}\right\} 
          \right\} & \text{ if } 0 \leq \ell_m < L_n ,  \\  
    \sigma\left\{\cG_{m-1}^{(n)} \bigcup \sigma 
           \left\{\eta_{k_m-i,L_n+i}\;|\;0 \leq i <k_m\right\} 
          \right\} & \text{ if } \ell_m = L_n . 
                 \end{cases} 
 \end{equation*}  
Furthermore, write
 ${M_m^{(n)}}:= (M_{m}^{1,(n)},M_{m}^{2,(n)},\ldots, M_{m}^{J+3,(n)})^{\top}$ 
 and 
 \begin{equation}  \label{osszevont} 
   \bar{\eta}_{k,L_n}:= \sum_{j=0}^{k-1} \eta_{j+1,k+L_n-j-1}, \qquad  
   \tilde{\eta}_{k,L_n}:= \sum_{j=0}^{k-1} \sum_{i=0}^{k+L_n-j-1} 
  \varrho_0^{k+L_n-j-i-1} \eta_{j+1,i}. 
 \end{equation}     
For the sake of convenience and better readability, in the following definition
 of $M_{m}^{i,(n)}$ we will simply write $\ell$ instead of $\ell_m$ and $k$
 instead of $k_m$.   
Define 
 \begin{equation*} 
  M_{m}^{1,(n)}
  := \begin{cases} 
    \frac{1}{n \beta_0}  
     \left[ \eta_{k,\ell}^2 - 1 + \eta_{k,\ell} 
        \left( \beta_0 \sum_{j=0}^{2\ell} \varrho_0^j  
               - \sum_{j=\ell}^J b_{0,j} \varrho_0^{j-\ell}  
        \right)  
     \right] 
       & \text{if } 0 \leq \ell < L_n , \\  
    \frac{1}{n \beta_0 k}             
      \left[ \bar{\eta}_{k,L_n}^2 - k 
      + \bar{\eta}_{k,L_n} 
       \left( \beta_0 \sum_{j=0}^{k-1} \sum_{i=0}^{2(k+L_n-j-1)} 
                \varrho_0^i  
       \right) 
      \right] 
       & \text{if } \ell=L_n , 
     \end{cases} 
 \end{equation*}
 \begin{equation*} 
  M_{m}^{2,(n)}
  := \begin{cases}
    \frac{1}{n}  
     \eta_{k,\ell} 
    \bigg[ \sum_{j=0}^{\ell-1} \varrho_0^{\ell-j-1} \eta_{k,j} 
        + \frac{\beta_0}{2} 
        \left( \sum_{j=0}^{\ell-1} \varrho_0^j \right)^2  
        + \frac{\beta_0}{2} \sum_{j=1}^{2\ell} j \varrho_0^{j-1} & \\   
    \qquad \qquad
  - \beta_0 \sum_{j=0}^J b_{0,j} \sum_{i=0}^{j \wedge (\ell-1)} 
              \varrho_0^{\ell+j-1-2i} 
  - \beta_0 \sum_{j=\ell+1}^J b_{0,j} (j-\ell) \varrho_0^{j-\ell-1} 
    \bigg]  & \\ 
        \text{\hskip8cm  if } 0 \leq \ell < L_n , & \\  
    \frac{1}{n k} \bar{\eta}_{k,L_n} 
      \bigg[ \tilde{\eta}_{k,L_n}  
            + \frac{\beta_0}{2} \sum_{j=0}^{k-1} 
              \left( \sum_{i=0}^{k+L_n-j-2} \varrho_0^i \right)^2 
            \! +  \! \frac{\beta_0}{2} \sum_{j=0}^{k-1} 
                 \sum_{i=1}^{2(k+L_n-j-1)} i\varrho_0^{i-1}  & \\  
    \qquad \qquad  
     -\sum_{i=0}^{k-1} \sum_{j=0}^J b_{0,j} 
       \sum_{n=0}^{j} \varrho_0^{j+k+L_n-i-1-2n}    
      \bigg] & \\ 
        \text{\hskip8cm if } \ell=L_n , & 
   \end{cases} 
 \end{equation*}
 and for $j = 0, \ldots, J$ 
 \begin{equation*} 
  M_{m}^{j+3,(n)}
  := \begin{cases}
    - \pmb{1}_{\{j \geq \ell\}} n^{-1/2} \varrho_0^{j-\ell} \eta_{k,\ell} 
       & \text{if } 0 \leq \ell < L_n , \\  
    0 
       & \text{if } \ell=L_n . 
   \end{cases} 
 \end{equation*}

Now, by the independence of the $\eta_{k,l}$'s, it is easy to see that 
 \[
   \EE \left( M_{m}^{(n)} | \cG_{m-1}^{(n)} \right) = 0
 \] 
 for $1 \leq m \leq K_n(L_n+1)$. 
Hence, we can see that $\left( M_{m}^{(n)} \right)_{ 1 \leq m \leq K_n(L_n+1) }$ are
 martingale differences with respect to the filtration
 $\left( \cG_{m}^{(n)} \right)_{ 0 \leq m \leq K_n(L_n+1) }$. 
Furthermore, recalling \eqref{kl}, \eqref{kl2}, \eqref{kL0} and \eqref{kL} we
 can see that for sufficiently large $n$ (s.t. $L_n > J$) we clearly have 
 \begin{equation*} 
  \mathbf{Z}_n = \sum_{m=1}^{K_n(L_n+1)} M_{m}^{(n)} . 
 \end{equation*} 

Next observe that the sequence consisting of the conditional covariances of
 $M_m^{(n)}$ 
 tends to $\Sigma$ in probability, i.e. 
 \begin{equation}  \label{kovi}   
    \gamma_{i,j}^{(n)} := 
    \sum_{m=1}^{K_n(L_n+1)} \EE \left(
     M_{m}^{i,(n)} M_{m}^{j,(n)} \;|\; 
     \cG_{m-1}^{(n)} \right) \stoch \sigma_{i,j} 
     \qquad \text{for } i,j = 1, 2, \ldots, J+3  
 \end{equation}
 as well as the conditional Liapounov condition holds, i.e.   
 \begin{equation}  \label{ljapunov} 
  \sum_{m=1}^{K_n(L_n+1)} \EE 
    \left( \left\|M_{m}^{(n)}\right\|^4 \;|\;  
     \cG_{m-1}^{(n)} \right) \stoch 0 . 
 \end{equation} 
In fact, we will show more: the convergence results in \eqref{kovi} and
 \eqref{ljapunov} are valid even in almost sure sense. 
For what follows (for the martingale limit theorem that we shall apply),
 however, the convergence in probability is sufficient. 

To show how to check \eqref{kovi} we only demonstrate two cases. 
Firstly, for $i=1$ and $j=2$ write 
 \begin{equation}  \label{de1}  
   \delta_{k,\ell} := \sum_{i=0}^{\ell-1} 
           \varrho_0^{\ell-i-1} \eta_{k,i} 
    +\frac{\beta_0}{2} \left( \sum_{i=1}^{\ell-1} \varrho_0^{i}
                       \right)^2 
    +\frac{\beta_0}{2} \sum_{i=1}^{2\ell} i \varrho_0^{i-1} , 
 \end{equation}  
 \begin{equation}  \label{de2} 
    \tilde{\delta}_{k,L_n} := 
       \tilde{\eta}_{k,L_n} 
       +\frac{\beta_0}{2} \sum_{j=0}^{k-1} \left( 
         \sum_{i=0}^{k+L_n-j-2} \varrho_0^i
                                           \right)^2 
       +\frac{\beta_0}{2} \sum_{j=0}^{k-1} 
             \sum_{i=1}^{2(k+L_n-j-1)} i \varrho_0^{i-1} . 
 \end{equation} 
Thus we obtain 
 \begin{equation*}  
  \begin{aligned} 
   \gamma_{1,2}^{(n)} 
   &= \frac{1}{\beta_0 n^2} \sum_{k=1}^{K_n} \sum_{\ell=0}^{L_n-1} 
       \delta_{k,\ell} \EE \eta_{k,\ell} \left( 
         \eta_{k,\ell}^2 - 1 
             + \beta_0 \eta_{k,\ell} \sum_{i=0}^{2\ell} \varrho_0^i
                                         \right) \\ 
   &\quad + \frac{1}{\beta_0 n^2} \sum_{k=1}^{K_n} \frac{1}{k^2} 
      \tilde{\delta}_{k,L_n} \EE \bar{\eta}_{k,L_n} 
                                 \left( 
        \bar{\eta}_{k,L_n}^2 - k + \beta_0 \bar{\eta}_{k,L_n} 
         \sum_{j=0}^{k-1} \sum_{i=0}^{2(k+L_n-j-1)} 
            \varrho_0^i 
                                 \right) + o(1) \\  
  &=\frac{KL}{\beta_0}         \left[ 
          \frac{\beta_0}{2\left(1-\varrho_0\right)^2} 
         +\frac{\beta_0}{2\left(1-\varrho_0\right)^2} 
                               \right] 
      \frac{\beta_0}{1-\varrho_0} \\ 
  &\quad +\frac{1}{\beta_0 n^2} \sum_{k=1}^{K_n} \frac{1}{k^2} 
      \left[ 
          \frac{\beta_0 k}{2\left(1-\varrho_0\right)^2} 
         +\frac{\beta_0 k}{2\left(1-\varrho_0\right)^2} 
      \right] \frac{\beta_0 k^2}{1-\varrho_0} + o(1) 
  \to \sigma_{1,2}   
  \end{aligned} 
 \end{equation*} 
 a.s.~as $n \to \infty$. 
Note that all the terms containing the market price of risk parameters vanish,
 i.e.~their order is $o(n^2)$, hence we omit to display these terms. 
Secondly, take $i,j \in \{0,1, \ldots , J\}$ and consider
 $\gamma_{i+3,\; j+3}^{(n)}$. 
Now we obtain 
 \[
   \gamma_{i+3,\; j+3}^{(n)} = n^{-1} K_n \sum_{\ell=0}^{i \wedge j} 
      \varrho_0^{i+j-2\ell} \EE \eta_{k,\ell}^2 
     \to \sigma_{i+3,j+3} \quad 
     \text{ a.s. as } n \to \infty .
 \] 
The remaining cases can be derived in a similar way. 

To show \eqref{ljapunov} notice that even 
 \[
   n^{z} \sum_{m=1}^{K_n(L_n+1)} \EE 
    \left( \left( M_{m}^{i,(n)}\right)^4 \;|\;  
     \cG_{m-1}^{(n)} \right)
 \] 
 has an almost sure limit, where $z=2$ for $i=1,2$, and $z=1$ for
 $3 \leq i \leq J+3$.  
This can be shown easily by the application of Lemmas \ref{lem1}, \ref{lem1_2},
 \ref{lem2} and their corollaries. 
From this \eqref{ljapunov} is immediate. 

Finally, it is known that according to the martingale limit theorem (see
 Theorem VIII.3.33.~in \cite{shir}) that \eqref{kovi} and \eqref{ljapunov} are
 together sufficient to imply (\ref{ll3}). 
\proofend 

\noindent {\bf Proof of Lemma \ref{anl1}.} 
First consider the case $i \vee j \leq 2$. 
Then one can easily show that 
 \[
   \frac{1}{n^2} \EE \partial_i \partial_j 
               \cL_n(\beta_0,\varrho_0,\mathbf{b}_0)  
         \to -\sigma_{i,j} 
  \quad \text{as } n \to \infty .
 \] 
Hence, by Lemma \ref{lem1}, Lemma \ref{lem2} (with $\kappa = 2$) and/or by the
 corollaries following from these lemmas one can easily see that 
 \[
   \frac{1}{n^2} \partial_i \partial_j 
               \cL_n(\beta_0,\varrho_0,\mathbf{b}_0) 
          \to -\sigma_{i,j} 
        \quad \text{a.s. as } n \to \infty .
 \] 
To demonstrate the method, consider the most complicated case, where $i=1$,
 $j=2$. 
Recalling notations \eqref{osszevont}, \eqref{de1} and \eqref{de2}  we have
 a.s. 
 \begin{equation*} 
  \begin{aligned} 
  \partial_1 \partial_2 
               \cL_n(\beta_0,\varrho_0,\mathbf{b}_0) 
   &= - \sum_{k=1}^{K_n} \sum_{\ell=0}^{L_n-1} \eta_{k,\ell} 
      \left[ \frac{2}{\beta_0} \delta_{k,\ell} 
             - \sum_{i=1}^{2\ell} i \varrho_0^{i-1} 
      \right] 
     + \sum_{k=1}^{K_n} \sum_{\ell=0}^{L_n-1} 
       \delta_{k,\ell} \sum_{i=1}^{2\ell} \varrho_0^{i} \\ 
   &\quad - \sum_{k=1}^{K_n} \frac{1}{k} \bar{\eta}_{k,L_n} 
       \left[ \frac{2}{\beta_0} \tilde{\delta}_{k,L_n} 
              - \sum_{j=0}^{k-1} 
             \sum_{i=1}^{2(k+L_n-j-1)} i \varrho_0^{i-1}
       \right]  \\ 
   &\quad + \sum_{k=1}^{K_n} \frac{1}{k} \tilde{\delta}_{k,L_n} 
         \sum_{j=0}^{k-1} 
             \sum_{i=0}^{2(k+L_n-j-1)} \varrho_0^{i} 
     + o(n^2) . 
  \end{aligned} 
 \end{equation*} 
Note that all the terms containing the market price of risk parameters vanish,
 i.e.~their order is $o(n^2)$, hence we omit to display these terms. 
For the expected values we have 
 \begin{equation*} 
  \begin{aligned}  
   \frac{1}{n^2} &\EE \partial_1 \partial_2 
               \cL_n(\beta_0,\varrho_0,\mathbf{b}_0)  
    = - \frac{\beta_0}{n^2}  \sum_{k=1}^{K_n} \sum_{\ell=0}^{L_n-1} 
         \left[ \sum_{i=0}^{2\ell} \varrho_0^i \right] 
                                      \left[ 
          \frac{1}{2} \left( \sum_{i=0}^{\ell-1} \varrho_0^i \right)^2 
         +\frac{1}{2} \sum_{i=0}^{2\ell} i \varrho_0^{i-1} 
                                      \right]  \\          
    &\quad - \frac{\beta_0}{n^2}  \sum_{k=1}^{K_n} \frac{1}{k} 
       \left[ \sum_{j=0}^{k-1} 
             \sum_{i=1}^{2(k+L_n-j-1)} \varrho_0^{i} 
       \right] \\  
    &\quad\quad \quad\quad \quad\quad \times 
        \left[ 
           \frac{1}{2} \sum_{j=0}^{k-1} \left( 
             \sum_{i=0}^{k+L_n-j-2} \varrho_0^i
                                        \right)^2 
           +\frac{1}{2} \sum_{j=0}^{k-1} 
             \sum_{i=1}^{2(k+L_n-j-1)} i \varrho_0^{i-1} 
         \right] + o(1)  \\ 
    &= - \frac{\beta_0 K L }{1-\varrho_0} 
        \left[ 
          \frac{1}{2} \frac{1}{\left(1-\varrho_0\right)^2} 
         +\frac{1}{2} \frac{1}{\left(1-\varrho_0\right)^2} 
        \right] \\ 
    &\quad -\frac{\beta_0 K^2}{2\left(1-\varrho_0\right)} 
        \left[ 
          \frac{1}{2} \frac{1}{\left(1-\varrho_0\right)^2} 
         +\frac{1}{2} \frac{1}{\left(1-\varrho_0\right)^2} 
        \right] + o(1) \to - \sigma_{1,2}  
  \end{aligned} 
 \end{equation*} 
 as $n \to \infty$. 

Now consider the case $i \vee j > 2$. 
Then in a similar way one can easily show that 
 \[
   \frac{1}{n} \EE \partial_i \partial_j 
               \cL_n(\beta_0,\varrho_0,\mathbf{b}_0)  
         \to \lambda  
  \quad \text{as } n \to \infty , 
 \] 
 where $\lambda \in \RR$ and $\lambda = -\sigma_{i,j}$ for $i \wedge j >2 $. 
Hence, by the application of Lemma \ref{lem1_2}, Lemma \ref{lem2} (with
 $\kappa = 1$) one can easily see that 
 \[
   \frac{1}{n} \partial_i \partial_j 
               \cL_n(\beta_0,\varrho_0,\mathbf{b}_0) 
          \to -\sigma_{i,j} 
        \quad \text{a.s. as } n \to \infty
 \] 
 for $i \wedge j > 2$, and 
 \[
   \frac{1}{n^{3/2}} \partial_i \partial_j 
               \cL_n(\beta_0,\varrho_0,\mathbf{b}_0) 
          \to 0 = -\sigma_{i,j} 
        \quad \text{a.s. as } n \to \infty 
 \] 
 for $i \wedge j \leq 2$. 
For instance, taking $i,j \in \{0,1,\ldots , J\}$ we have 
 \[
   \partial_{i+3} \partial_{j+3}  
               \cL_n(\beta_0,\varrho_0,\mathbf{b}_0) 
    = - K_n \sum_{\ell=0}^{i \wedge j} \varrho_0^{i+j-2\ell} ,
 \] 
 from which the statement is immediate. 
The remaining cases can easily be calculated in a similar way. 
\proofend

\noindent {\bf Proof of Lemma \ref{anl2}.} 
Based on Lemma \ref{lem1}, Lemma \ref{lem1_2}, Lemma \ref{lem2} one can show
 that 
 \[
   \frac{1}{2} n^{l_{i,j_2}}\partial_i \partial_{j_1} 
     \partial_{j_2} \cL_n \left( \beta, \varrho, 
     \mathbf{b} \right)
 \] 
 has an almost sure limit uniformly in $( \beta, \varrho, \mathbf{b} ) \in H$. 
This can be shown similarly to the uniform convergence in \eqref{fosor}. 
(Recall also Remark \ref{uniform} for this. 
 Notice that in fact the higher order derivatives of the likelihood function
  will have at most the same speed of convergence as the first order ones in
  their asymptotic expansion due to its relatively easy dependence on the
  parameters.)
Thus, by Lemma \ref{inverz} we obtain that $r^{n,i}_{j_1,j_2}$ is stochastically
 bounded for all $1 \leq j_1, j_2 \leq J+3$. 
On the other hand, the estimators $\hbeta_n, \hvarrho_n, \hbb_n$ are proved to
 be strongly consistent and thus $D_n^{\top} \to 0$ a.s.~as
 $n \to \infty$. 
Hence we obtain that $B_n^i = D_n^{\top} R_n^i$ converges to zero in probability. 
\proofend

\section*{Appendix B}  
\renewcommand{\theequation}{B.\arabic{equation}}
\renewcommand{\thesection}{B}
\renewcommand{\thetheorem}{B.\arabic{theorem}}
\renewcommand{\thelemma}{B.\arabic{lemma}}
\renewcommand{\theproposition}{B.\arabic{proposition}}
\renewcommand{\thedefinition}{B.\arabic{definition}}
\renewcommand{\thecorollary}{B.\arabic{corollary}}
\renewcommand{\theremark}{B.\arabic{remark}}

\setcounter{equation}{0}
\setcounter{theorem}{0}
\setcounter{lemma}{0}
\setcounter{proposition}{0}
\setcounter{definition}{0}
\setcounter{remark}{0}
\setcounter{corollary}{0}

In what follows we summarise some simple but useful lemmas that are often used
 in the proofs of the main results. 
They give some general statements which are not model specific (which 
 was the reason for presenting them in a separate appendix). 

\begin{lemma}  \label{lem1} 
Let $\xi_{k,\ell,n}$, $k, \ell, n \in \NN$, be random variables such that for
 each $n \in \NN$ the sets $\{ \xi_{k,\ell,n} : \ell \in \NN \}$, $k \in \NN$,
 are independent (i.e., the $\sigma$-algebras
 $\sigma( \xi_{k,\ell,n} : \ell \in \NN )$, $k \in \NN$, are independent), and
 $\sup_{ k, \ell, n \in \NN } \EE \xi_{k,\ell,n}^8 < \infty$. 
Let $K_n, L_n$, $n \in \NN$, be positive integers such that $K_n = nK + o(n)$
 and $L_n = n L + o(n)$ as $n \to \infty$ with some $K > 0$ and $L > 0$.  
Then 
 \[
   n^{-2} \sum_{k=1}^{K_n} \sum_{\ell=1}^{L_n}
   \left( \xi_{k,\ell,n}^{2} - \EE \xi_{k,\ell,n}^{2} \right)
   \to 0 
   \qquad \text{a.s.~as $n \to \infty$.} 
 \]
\end{lemma} 

\noindent
{\bf Proof.} 
It suffices to show that for all $\vare > 0$ we have
 \[
   \sum_{n=1}^\infty \PP( | \zeta_n | > \vare n^2 ) < \infty ,
 \]
 where
 \[
   \zeta_n := \sum_{k=1}^{K_n} \sum_{\ell=1}^{L_n}
               \left( \xi_{k,\ell,n}^2 - \EE \xi_{k,\ell,n}^2 \right) .
 \]
By Markov inequality we obtain
 $\PP( | \zeta_n | > \vare n^2 ) \leq \vare^{-4} n^{-8} \EE \zeta_n^4$, hence it
 is enough to show that $\EE \zeta_n^4 = O(n^{7-\delta})$ as $n\to\infty$ with
 some $\delta>0$.
We have
 \[
   \EE \zeta_n^4 = \sum_{k_1,\,k_2,\,k_3,\,k_4\,=\,1}^{K_n} \;
                    \sum_{\ell_1,\,\ell_2,\,\ell_3,\,\ell_4\,=\,1}^{L_n}
                     \EE \zeta_{k_1,\ell_1,n} \zeta_{k_2,\ell_2,n}
                         \zeta_{k_3,\ell_3,n} \zeta_{k_4,\ell_4,n} ,
 \]
 where $\zeta_{k,\ell,n} := \xi_{k,\ell,n}^2 - \EE \xi_{k,\ell,n}^2$.
By the Cauchy--Schwartz inequality
 \[
   | \EE \zeta_{k_1,\ell_1,n} \zeta_{k_2,\ell_2,n}
         \zeta_{k_3,\ell_3,n} \zeta_{k_4,\ell_4,n} |
   \leq \big( \EE \zeta_{k_1,\ell_1,n}^4 \EE \zeta_{k_2,\ell_2,n}^4
              \EE \zeta_{k_3,\ell_3,n}^4 \EE \zeta_{k_4,\ell_4,n}^4 \big)^{1/4} .
 \] 
Moreover
 \[
   \EE \zeta_{k,\ell,n}^4
   =\EE ( \xi_{k,\ell,n}^2 - \EE \xi_{k,\ell,n}^2 )^4
   \leq 2^3 \left( \EE \xi_{k,\ell,n}^8 + ( \EE \xi_{k,\ell,n}^2 )^4 \right)
   \leq 16 \EE \xi_{k,\ell,n}^8 
   \leq 16 M_8 ,
 \]
 where $M_8 := \sup_{k, \ell, n \in \NN} \EE \xi_{k,\ell,n}^8 < \infty$ by the
 assumptions.
Hence we conclude
 \[
   | \EE \zeta_{k_1,\ell_1,n} \zeta_{k_2,\ell_2,n}
         \zeta_{k_3,\ell_3,n} \zeta_{k_4,\ell_4,n} |
   \leq 16 M_8.
 \]
By the assumptions the sets $\{ \zeta_{k,\ell,n} : \ell \in \NN \}$, $k \in \NN$,
 are independent for each $n \in \NN$, and $\EE \zeta_{k,\ell,n} = 0$ for all 
 $k, \ell, n \in \NN$, hence
 \begin{equation} 
  \begin{aligned} 
   \EE \zeta_n^4
   = \sum_{\ell_1,\,\ell_2,\,\ell_3,\,\ell_4\,=\,1}^{L_n}
     \Bigg( &\sum_{k=1}^{K_n} \EE \zeta_{k,\ell_1,n} \zeta_{k,\ell_2,n}
                                 \zeta_{k,\ell_3,n} \zeta_{k,\ell_4,n}  \\ 
            & + 6 \sum_{1 \leq k_1 < k_2 \leq K_n}
                \EE \zeta_{k_1,\ell_1,n} \zeta_{k_1,\ell_2,n}
                    \zeta_{k_2,\ell_3,n} \zeta_{k_2,\ell_4,n} \Bigg) .
  \end{aligned} 
 \end{equation} 
Consequently we obtain $\EE \zeta_n^4 = O(n^6)$ as $n\to\infty$. 
\proofend

\begin{lemma}  \label{lem1_2} 
Let $\xi_{k,n}$, $k, n \in \NN$, be random variables such that for each
 $n \in \NN$ the random variables $\xi_{k,n}$, $k \in \NN$, are independent and
 $\sup_{k, n \in \NN} \EE \xi_{k,n}^{4} < \infty$. 
Let $K_n$, $n \in \NN$, be positive integers such that $K_n = nK + o(n)$ as
 $n \to \infty$ with some $K > 0$. 
Then 
 \[
   n^{-1} \sum_{k=1}^{K_n} 
          \left( \xi_{k,n} - \EE \xi_{k,n} \right)
   \to 0 
   \qquad \text{a.s.~as $n \to \infty$.} 
 \]
\end{lemma} 

\noindent
{\bf Proof.} \  
This statement can be proved almost readily in the same way as Lemma
 \ref{lem1}. 
\proofend

\begin{corollary}  \label{cor1} 
Let $\xi_{k,\ell,n}$ and $\zeta_{k,\ell,n}$, $k, \ell, n \in \NN$, be random
 variables such that for each $n \in \NN$ the sets 
 $\{ \xi_{k,\ell,n}, \zeta_{k,\ell,n} : \ell \in \NN \}$, $k \in \NN$, are
 independent (i.e., the $\sigma$-algebras
 $\sigma( \xi_{k,\ell,n}, \zeta_{k,\ell,n} : \ell \in \NN )$, $k \in \NN$, are
 independent), and  
 $\sup_{k, \ell,n \in \NN} \EE ( \xi_{k,\ell,n}^8 + \zeta_{k,\ell,n}^8 ) < \infty$. 
Let $K_n, L_n$, $n \in \NN$, be positive integers such that $K_n = nK + o(n)$
 and $L_n = n L + o(n)$ as $n \to \infty$ with some $K > 0$ and $L > 0$.
Then 
 \[
   n^{-2} \sum_{k=1}^{K_n} \sum_{\ell=1}^{L_n}
    ( \xi_{k,\ell,n} \zeta_{k,\ell,n} - \EE \xi_{k,\ell,n} \zeta_{k,\ell,n} )
   \to 0 
   \qquad \text{a.s.~as $n \to \infty$.} 
 \]
\end{corollary} 

\noindent
{\bf Proof.} \
Clearly
 \begin{align*}
  &\xi_{k,\ell,n} \zeta_{k,\ell,n} - \EE \xi_{k,\ell,n} \zeta_{k,\ell,n} \\
  &= \frac{1}{4}
    \bigg[ \Big\{ ( \xi_{k,\ell,n} + \zeta_{k,\ell,n})^2
                    - \EE ( \xi_{k,\ell,n} + \zeta_{k,\ell,n})^2 \Big\}
          - \Big\{ ( \xi_{k,\ell,n} - \zeta_{k,\ell,n})^2
                   - \EE ( \xi_{k,\ell,n} - \zeta_{k,\ell,n})^2 \Big\} \bigg] ,
 \end{align*}
 and we can apply Lemma \ref{lem1} for 
 $\{ \xi_{k,\ell,n} + \zeta_{k,\ell,n} : k, \ell, n \in \NN \}$
 and $\{ \xi_{k,\ell,n} - \zeta_{k,\ell,n} : k, \ell, n \in \NN \}$.
\proofend

\begin{corollary}  \label{cor2} 
Let $\xi_{k,\ell,n}$, $k, \ell, n \in \NN$, be random variables such that for
 each $n \in \NN$ the sets $\{ \xi_{k,\ell,n} : \ell \in \NN \}$, $k \in \NN$,  
 are independent (i.e., the $\sigma$-algebras
 $\sigma( \xi_{k,\ell,n} : \ell \in \NN )$, $k \in \NN$, are independent), and 
 $\sup_{k, \ell,n \in \NN} \EE \xi_{k,\ell,n}^8 < \infty$. 
Let $K_n, L_n$, $n \in \NN$, be positive integers such that $K_n = nK + o(n)$
 and $L_n = n L + o(n)$ as $n \to \infty$ with some $K > 0$ and $L > 0$. 
Then 
 \[
   n^{-2} \sum_{k=1}^{K_n} \sum_{\ell=1}^{L_n}
    ( \xi_{k,\ell,n} - \EE \xi_{k,\ell,n} )
   \to 0 
   \qquad \text{a.s. as \  $n \to \infty$.} 
 \]
\end{corollary} 

\noindent
{\bf Proof.} \
Corollary \ref{cor1} applies with $\zeta_{k,\ell,n} = 1$, $k, \ell, n \in\NN$.
\proofend

\begin{lemma}  \label{lem2} 
Let $\kappa \in \{ 1, 2 \}$. 
Let $\xi_{k,j,n}$, $k, j, n \in \NN$, be random variables such that for each
 $n \in \NN$ the sets $\{ \xi_{k,j,n} : k \in \NN \}$, $j \in \NN$, are
 independent (i.e., the $\sigma$-algebras $\sigma( \xi_{k,j,n} : k \in \NN )$,
 $j \in \NN$, are independent), and
 $\sup_{k, j, n \in \NN} \EE \xi_{k,j,n}^{4\kappa} < \infty$. 
Let $K_n$, $n \in \NN$, be positive integers such that $K_n = nK + o(n)$ as
 $n \to \infty$ with some $K > 0$. 
Then 
 \[
   n^{-\kappa} \sum_{k=1}^{K_n} k^{-1}
    \left[ \left( \sum_{j=1}^k \xi_{k,j,n} \right)^{\kappa} 
           - \EE \left( \sum_{j=1}^k \xi_{k,j,n} \right)^{\kappa} \right]
   \to 0 
   \qquad \text{a.s.~as $n \to \infty$.} 
 \]
\end{lemma} 

\noindent
{\bf Proof.} \  
Consider the case $\kappa = 2$. 
Clearly
 \[
   \left( \sum_{j=1}^k \xi_{k,j,n} \right)^2
          - \EE \left( \sum_{j=1}^k \xi_{k,j,n} \right)^2
   = 
   \sum_{j_1=1}^k \sum_{j_2=1}^k 
    ( \xi_{k,j_1,n} \xi_{k,j_2,n} - \EE \xi_{k,j_1,n} \xi_{k,j_2,n} ) .
 \]
As in the proof of Lemma \ref{lem1} it suffices to show that
 $\EE \zeta_n^4 = O ( n^{7-\delta})$ as $n \to \infty$ with some $\delta > 0$,
 where
 \[
   \zeta_n := \sum_{k=1}^{K_n} k^{-1} \sum_{j_1=1}^k \sum_{j_2=1}^k \zeta_{k,j_1,j_2,n}
 \]
 with
 \[
   \zeta_{k,j_1,j_2,n} := \xi_{k,j_1,n} \xi_{k,j_2,n} - \EE \xi_{k,j_1,n} \xi_{k,j_2,n} .
 \]
We have
 \begin{equation} 
  \begin{aligned} 
   \zeta_n^4
   = &\sum_{k_1,\,k_2,\,k_3,\,k_4\,=\,1}^{K_n} \Bigg[ 
             (k_1k_2k_3k_4)^{-1} \\ 
     &\qquad \times 
      \sum_{j_1,\,j_2\,=\,1}^{k_1} \; \sum_{j_3,\,j_4\,=\,1}^{k_2} \;
      \sum_{j_5,\,j_6\,=\,1}^{k_3} \; \sum_{j_7,\,j_8\,=\,1}^{k_4} 
       \EE \zeta_{k_1,j_1,j_2,n} \zeta_{k_2,j_3,j_4,n}
           \zeta_{k_3,j_5,j_6,n} \zeta_{k_4,j_7,j_8,n} \Bigg] .
  \end{aligned} 
 \end{equation} 
As in the proof of Lemma \ref{lem1} we obtain 
 \[ 
   | \EE \zeta_{k_1,j_1,j_2,n} \zeta_{k_2,j_3,j_4,n}
         \zeta_{k_3,j_5,j_6,n} \zeta_{k_4,j_7,j_8,n} |
   \leq 16 M_8,
 \]
 where $M_8 := \sup_{k, j, n \in \NN} \EE \xi_{k,j,n}^8 < \infty$ by the assumptions. 
By the independence of the sets $\{ \xi_{k,j,n} : k \in \NN \}$, $j \in \NN$, we
 obtain that
 $\EE \zeta_{k_1,j_1,j_2,n} \zeta_{k_2,j_3,j_4,n} \zeta_{k_3,j_5,j_6,n} \zeta_{k_4,j_7,j_8,n}
   = 0$ 
 if one of the sets $\{j_1,j_2\}$, $\{j_3,j_4\}$, $\{j_5,j_6\}$, $\{j_7,j_8\}$ is
 disjoint from the other three sets.
Consequently
 \[
   \sum_{j_1,\,j_2\,=\,1}^{k_1} \; \sum_{j_3,\,j_4\,=\,1}^{k_2} \;
    \sum_{j_5,\,j_6\,=\,1}^{k_3} \; \sum_{j_7,\,j_8\,=\,1}^{k_4} 
     \EE \zeta_{k_1,j_1,j_2,n} \zeta_{k_2,j_3,j_4,n}
         \zeta_{k_3,j_5,j_6,n} \zeta_{k_4,j_7,j_8,n}
   = O ( n^6 )
 \]
 as $n \to \infty$.
Using $\sum_{k=1}^n k^{-1} = O (\log n)$ we conclude
 $\EE \zeta_n^4 = O \big( n^6 (\log n)^4 \big)$ as $n \to \infty$.

The case $\kappa = 1$ can be proved almost readily in the same way. 
\proofend

\begin{corollary}  \label{cor3} 
Let $\xi_{k,j,n}$, $\zeta_{k,j,n}$, $k, j, n \in \NN$, be random variables such
 that for each $n \in \NN$ the sets $\{ \xi_{k,j,n} : k \in \NN \}$,
 $j \in \NN$, are independent (i.e., the $\sigma$-algebras 
 $\sigma( \xi_{k,j,n}, \zeta_{k,j,n} : k \in \NN )$, $j \in \NN$, are
 independent), and
 $\sup_{k, \ell,n \in \NN} \EE ( \xi_{k,\ell,n}^8 + \zeta_{k,\ell,n}^8 ) < \infty$. 
Let $K_n$, $n \in \NN$, be positive integers such that
 $K_n = nK + o(n)$ as $n \to \infty$ with some $K > 0$. 
Then 
 \[
   n^{-2} \sum_{k=1}^{K_n} k^{-1}
   \sum_{j=1}^k
    ( \xi_{k,j,n} \zeta_{k,j,n} - \EE \xi_{k,j,n} \zeta_{k,j,n} ) 
   \to 0 
   \qquad \text{a.s. as \  $n \to \infty$.} 
 \]
\end{corollary} 

\noindent
{\bf Proof.} \
Similar to the proof of Corollary \ref{cor1}.
\proofend

\begin{lemma}  \label{inverz} 
Let $H_n: \RR^{n+1} \to \RR$, $n \in \NN$, be measurable functions and 
  $\{ \xi_n \}_{n \in \NN}$ be a sequence of random variables. 
Suppose that $H: \RR \to \RR$ is continuous and $C$ is a compact subset of
 $\RR$ such that 
 \[
   \sup_{\alpha \in C} 
    \left| H_n(\xi_1, \ldots , \xi_n, \alpha) - H(\alpha) \right| \to 0  
   \qquad \text{$\PP$-a.s.}  
 \] 
Then, given random variables $\alpha_n$, $n \in \NN$, with
 $\PP\left( \alpha_n \in C \right) = 1$, the sequence 
 \[
   \{ H_n(\xi_1, \ldots , \xi_n, \alpha_n) \}_{n \in \NN}
 \] 
 is stochastically bounded in the following sense: 
 \begin{equation}  \label{sbdef}
   \lim_{R \to \infty} \limsup_{n \to \infty} 
       \PP \left( | H_n (\xi_1, \ldots , \xi_n, \alpha_n) | 
                   > R 
           \right) = 0. 
 \end{equation} 
\end{lemma}
 
{\bf Proof.} 
Let $R>0$. 
We have 
 \begin{equation}  \label{sbfelb} 
  \begin{aligned}
   \PP  ( | H_n (&\xi_1, \ldots , \xi_n, \alpha_n) | > R 
        )  \\ 
   &\leq \PP \left( | H_n (\xi_1, \ldots , \xi_n, \alpha_n) 
              - H (\alpha_n) | + | H (\alpha_n) | > R \right) \\ 
   &\leq \PP \left( | H_n (\xi_1, \ldots , \xi_n, \alpha_n) 
              - H (\alpha_n) | > \frac{R}{2} \right) 
       + \PP \left( | H (\alpha_n) | > \frac{R}{2} \right)  \\ 
   &\leq \PP \left( \sup_{\alpha \in C} 
                      | H_n (\xi_1, \ldots , \xi_n, \alpha) 
              - H (\alpha) | > \frac{R}{2} \right) 
       + \PP \left( \sup_{\alpha \in C} | H (\alpha) | > \frac{R}{2} \right), 
  \end{aligned} 
 \end{equation}  
 hence, taking in both sides of (\ref{sbfelb}) first the `$\limsup$' as
 $n \to \infty$ and then the limit as $R \to \infty$, one gets the desired
 statement. 
\proofend 

Let us remark that stochastic boundedness is not necessarily defined as in
 \eqref{sbdef} in the literature. 
However, it suffices for our purpose. 
For this, we note that given a sequence of random variables, say
 $\{X_n\}_{n \in \NN}$, with limit $0$ in the sense of convergence in probability
 and given a stochastically bounded (in the sense of (\ref{sbdef}))
 sequence of r.v.'s, say $\{Y_n\}_{n \in \NN}$, one easily gets that 
 $ X_n \cdot Y_n $ converges to $0$ in probability, as well. 

\begin{lemma}  \label{anl4} 
Let $\{X_n\}_{n \in \ZZ_+}$, $\{Y_n\}_{n \in \ZZ_+}$ and $\{Z_n\}_{n \in \ZZ_+}$ be
 sequences of random matrices of type $m\times1$, $m\times m$ and $m\times1$,
 respectively, such that $X_n=Y_nZ_n$, where $m \in \NN$. 
Suppose that $X_n \distr X$, $Y_n \stoch A$, where $A$ is a non-degenerate
 matrix. 
Then $Z_n \distr A^{-1} X$. 
\end{lemma}
 
{\bf Proof.} 
Define 
 \begin{equation*} 
  Y_n^{\circleddash} := \begin{cases} 
     Y_n^{-1} & \text{if $Y_n$ invertible,} \\   
     0        & \text{otherwise.}  
                     \end{cases} 
 \end{equation*} 
Clearly, we have $\PP(Y_n^{\circleddash}=Y_n^{-1}) \to 1$, as $n \to \infty$. 
Hence, 
 \[
   \PP(||Y_n^{\circleddash}-A^{-1}|| \geq \vare) 
    \leq \PP(||Y_n^{-1}-A^{-1}|| \geq \vare \text{ and } 
               Y_n^{\circleddash}=Y_n^{-1}) 
       + \PP(Y_n^{\circleddash} \neq Y_n^{-1}) 
     \to 0
 \] 
 as $n \to \infty$, that is $Y_n^{\circleddash} \stoch A^{-1}$. 
By Slutsky's Lemma, $Y_n^{\circleddash} X_n \stoch A^{-1} X$. 
Further, $\PP(Z_n = Y_n^{\circleddash} X_n) \to 1$, and hence $Z_n \stoch A^{-1} X$. 
Thus, the proof of Lemma \ref{anl4} is completed, and so is the proof of
 Theorem \ref{stabil}. 
\proofend 


\begin{thebibliography}{99}


\bibitem{FP} {\sc F\"ul\"op, E.} and {\sc Pap, G.} (2007),  
{\em Asymptotically optimal tests for a discrete time random field HJM
 type interest rate model},     
 Acta Scientiarum Mathematicarum, {\bf 73(3-4)}, 637--661.

\bibitem{FP2} {\sc F\"ul\"op, E.} and {\sc Pap, G.} (2009),  
{\em Strong consistency of maximum likelihood estimators for 
 a discrete-time random field HJM-type interest rate model},  
 Lithuanian Math. J., {\bf 49(1)}, 5-25. 

\bibitem{FE} {\sc F\"ul\"op, E.} (2009),  
{\em Simulations of a discrete time HJM type forward interest rate model},  
 unpublished manuscript. 

\bibitem{GPZ3} {\sc G\'all, J.}, {\sc Pap, G.} and 
  {\sc Zuijlen, M. v.} (2004), 
 {\em Maximum likelihood estimator of the volatility of forward rates driven 
  by geometric spatial AR sheet},   
 Journal of Applied Mathematics {\bf 2004(4)}, 293--309.

\bibitem{GPZ} {\sc G\'all, J.}, {\sc Pap, G.} and {\sc Zuijlen, M. v.}
 (2006), 
 {\em Forward interest rate curves in discrete time settings driven 
  by random fields}, 
 Computers \& Mathematics with Applications, {\bf 51(3-4)}, 387--396.  

\bibitem{GPP} {\sc G\'all, J.}, {\sc Pap, G.} and {\sc Peeters, W.} 
 (2007),  
{\em Random field forward interest rate models,
 market price of risk and their statistics}, 
 Analli dell'Universita di Ferrara Sez.~VII Sci.~Mat., 
 {\bf 53}, 233--242.

\bibitem{hm} {\sc Heijmans, R. D. H.} and {\sc Magnus, J. R.} 
    (1986), 
 {\em Consistent maximum-likelihood estimation with 
      dependent observations, The general (non-normal) 
      case and the normal case}, 
 Journal of Econometrics, {\bf 32}, 253--285. 

\bibitem{hm2} {\sc Heijmans, R. D. H.} and {\sc Magnus, J. R.} 
    (1986), 
 {\em Asymptotic normality of maximum likelihood estimators 
      obtained from normally distributed but dependent 
      observations}, 
 Econometric Theory, {\bf 2}, 374--412. 

\bibitem{shir} {\sc Jacod, J.}, {\sc Shiryayev, A. N.} (1987),  
 {\em Limit Theorems for Stochastic Processes}, Springer-Verlag, 
 Berlin. 

\bibitem{L} {\sc Le Cam, L.} (1986), 
 {\em Asymptotic Methods in Statistical Decision Theory}, 
 New York: Springer-Verlag. 

\bibitem{WP} {\sc Peeters, W.} (2008),  
 {\em Volatility estimation for different structures of random field
 interest rate models in discrete time},  
 Publicationes Mathematicae Debrecen, 
 {\bf 72(3-4)}, 317--334.

\bibitem{vdV} {\sc Vaart, A. W. van der} (1998), 
 {\em Asymptotic Statistics}, 
 Cambridge University Press. 

\bibitem{Y1} {\sc Ying, Z.} (1993), 
 {\em Asymptotic properties of a maximum likelihood estimator with data from a
  Gaussian process}, 
 Journal of Multivariate Analysis, {\bf 36}, 280--296.

\bibitem{Y2} {\sc Ying, Z.} (1993), 
 {\em Maximum likelihood of parameters under spatial sampling scheme}, 
 The Annals of Statistics, {\bf 21}, 1567--1590.

\end{thebibliography}
\end{document}